\documentclass[12pt]{elsarticle}

\newtheorem{thm}{Theorem} % [section]
 \newtheorem{cor}[thm]{Corollary}
 \newtheorem{lem}[thm]{Lemma}

 \newdefinition{rem}{Remark}
 
\newproof{proof}{Proof}
\newtheorem{defn}{Definition}

\usepackage{amsmath}
\usepackage{amsfonts}
\usepackage{amssymb}

\newcommand{\abs}[1]{\ensuremath{\left| #1 \right| }}
\newcommand{\qPs}{$q$-Pochhammer symbol}

\newcommand{\lhs}{left hand side}

\newcommand{\wrt}{with respect to}

\newcommand{\RRis}{Rogers--Ramanujan identities}

\newcommand{\ci}{cocycle identity}

\date{}
\journal{arxiv}

\begin{document}

\begin{frontmatter}

\title%[Multiple Combinatorial Numbers]
{Multiple Bracket Function, Stirling Number, and Lah Number Identities} 

\author{Hasan Coskun}
\ead{hasan.coskun@tamuc.edu}
%\url{http://faculty.tamuc.edu/hcoskun}
\address{Department of Mathematics, Texas A\&M
  University--Commerce, Binnion Hall, Room 314, Commerce, TX 75429}  
%\curraddr{Department of Mathematics, Texas A\&M
%  University--Commerce, Binnion Hall, Room 314, Commerce, TX 75429}  
%\email{hasan.coskun@tamu-commerce.edu}

\begin{abstract}
The author has constructed multiple analogues of several families of combinatorial numbers in a recent article, including the bracket symbol, and the Stirling numbers of the first and second kind. In the present paper, a multiple analogue of another sequence, the Lah numbers, is developed, and certain associated identities and significant properties of all these sequences are constructed. 
\end{abstract}

\begin{keyword}
multiple combinatorial numbers \sep multiple binomial coefficients \sep multiple bracket function \sep multiple Stirling numbers \sep multiple Lah numbers \sep well poised rational Macdonald functions 

\MSC 05A10
\sep 11B65
\sep 33D67
\end{keyword}

\end{frontmatter}

%\maketitle

\section{Introduction}
\label{introduction}
The Stirling numbers of the first and the second kind have been studied and their properties have been investigated extensively in number theory, combinatorics and other areas. 
One dimensional generalizations of these numbers have also been subject of interest. An important class of generalizations % of special numbers 
is their one parameter $q$-extensions. 
Many have made significant contributions to the $q$-extensions investigating their properties and applications.
We will give references to some of these important works in Section~\ref{mulipleSpecial} below.

In a recent paper, the author took a major step from one dimensional generalizations to multiple analogues of combinatorial numbers by constructing elegant multiple $qt$-generalizations of Stirling numbers of the first and second kind, besides sequences of other combinatorial numbers including multiple binomial, Fibonacci, Bernoulli, Catalan, and Bell numbers~\cite{CoskunDiscM}. In this paper, we focus on the multiple analogies of the factorial function, Stirling numbers of both kinds, and Lah numbers, and give interesting new identities they satisfy.

The multiple generalizations developed in~\cite{CoskunDiscM} are given in terms of the $qt$-binomial coefficients constructed in the same paper. Its definition may be written in the general form as
\begin{equation*}
%\label{qtbinom}
\binom{z}{\mu}_{\!\!\!q,t} := \dfrac{ q^{|\mu|} t^{2n(\mu)+(1-n)|\mu| } }  { (qt^{n-1} )_\mu} \prod_{1\leq i<j \leq n} \left\{\dfrac{ (qt^{j-i})_{\mu_i-\mu_j} } {(qt^{j-i-1})_{\mu_i-\mu_j}  } \right\}  %W^{\nearrow \hspace*{-6pt} s}
w_\mu(q^z t^{\delta(n)}; q, t)
\end{equation*}
where $\mu$ is a partition of at most $n$ parts, $z\in\mathbb{C}^n$ and $q,t\in\mathbb{C}$. The $w_\mu$ %$W^{\nearrow \hspace*{-6pt} s}_\mu$ 
function that enters the definition is a limiting case of the $BC_n$ well--poised symmetric rational Macdonald function $W_\lambda$. 
Note that this definition makes sense even when $\mu$ is not an integer partition, but is a vector $\mu\in \mathbb{C}^n$. The $qt$-binomial coefficients are constructed independently by Kaneko~\cite{Kaneko1} in a special case, and Okounkov~\cite{Okounkov1} using difference operator methods for integer partitions. Lassalle constructed an equivalent set of multiple binomial coefficient independently~\cite{Lassalle1}, and Sahi has developed the non-symmetric version~\cite{Sahi1}. The most general definition~(\ref{qtbinom}) of multiple $qt$-binomial coefficients we constructed in~\cite{CoskunDiscM} where $z\in\mathbb{C}^n$ (and $\mu\in\mathbb{C}^n$) are variables is essential in the  development of the multiple brackets, and thus the Stirling and Lah numbers in this paper.

\section{Background}
\label{back}

The (basic) \qPs\ $(a;q)_\alpha$ may be defined formally for complex parameters $q, \alpha\in\mathbb{C}$ as
\begin{equation}
\label{qPochSymbol}
(a)_\alpha = (a;q)_\alpha :=\dfrac{(a;q)_\infty}{(aq^\alpha;q)_\infty}
\end{equation} 
where the infinite product $(a;q)_\infty$ is defined by $(a;q)_\infty:=\prod_{i=0}^{\infty} (1-aq^i)$. Note that when $\alpha=m$ is a positive integer, the definition reduces to the finite product $(a;q)_m= \prod_{k=0}^{m-1}(1-aq^k)$. An elliptic
analogue 
is defined~\cite{FrenkelT1} by
\begin{equation}
  (a; q,p)_m := \prod_{k=0}^{m-1} \theta(aq^m)
\end{equation}
where $a\in \mathbb{C}$, $m$ is a positive integer, and the normalized theta
function $\theta(x)$ is given by
\begin{equation}
  \theta(x) = \theta(x;p) := (x; p)_\infty (p/x; p)_\infty
\end{equation}
for $x, p\in \mathbb{C}$ with $\abs{p}<1$. The definition is extended to negative $m$ by setting $(a; q,p)_m = 1/ (aq^{m}; q, p)_{-m} $. It is clear that when $p=0$, the elliptic $(a; q,p)_m$ reduces to the basic (trigonometric) \qPs~(\ref{qPochSymbol}).

For any partition $\lambda = (\lambda_1, \ldots, \lambda_n)$ and
$t\in\mathbb{C}$, define~\cite{Warnaar2}
\begin{equation}
\label{ellipticQtPocSymbol}
  (a)_\lambda=(a; q, p, t)_\lambda := \prod_{i=1}^{n} (at^{1-i};
  q,p)_{\lambda_i} .
\end{equation}
Note that when $\lambda=(\lambda_1) = \lambda_1$ is a single part
partition, then $(a; q, p, t)_\lambda = (a; q, p)_{\lambda_1} =
(a)_{\lambda_1}$. For brevity of notation, we also use
\begin{equation}
  (a_1, \ldots, a_k)_\lambda = (a_1, \ldots, a_k; q, p, t)_\lambda :=
  (a_1)_\lambda \ldots (a_k)_\lambda .
\end{equation}
Recall that we use $V$ to denote~\cite{Coskun1} the 
space of infinite lower--triangular matrices whose entries are rational
functions over the field
$\mathbb{F}=\mathbb{C}(q, p, t,r,a,b)$ which are indexed by partitions with respect to the partial inclusion
ordering $\subseteq$ defined by
\begin{equation}
\label{partialordering}
\mu \subseteq \lambda \;\Leftrightarrow \;\mu_i \leq \lambda_i, \quad
\forall i\geq 1.
\end{equation} 
The condition that a matrix $u\in V$ is lower triangular \wrt\ $\subseteq$ 
can be stated in the form
\begin{equation}
  u_{\lambda\mu} = 0 ,\, \quad \mathrm{when}\; \mu \not \subseteq
  \lambda.
\end{equation}
The multiplication operation defined by 
\begin{equation}
\label{multiplication}
  (uv)_{\lambda\mu} := \sum_{\mu\subseteq\nu\subseteq\lambda}
  u_{\lambda\nu} v_{\nu\mu}
\end{equation}
for matrices $u,v\in V$ makes $V$ into an $\mathbb{F}$-algebra. 

\subsection{Well--poised Macdonald functions}

The construction of our multiple Stirling numbers involves the elliptic
well--poised Macdonald functions $W_{\lambda/\mu}$
on $BC_n$~\cite{Coskun1}. 
These remarkable families of symmetric rational functions are first introduced in the author's Ph.D. thesis~\cite{Coskun0} in the basic (trigonometric) case, and later in~\cite{CoskunG1} in the more general elliptic form.   

Let $\lambda=(\lambda_1, \ldots, 
\lambda_n)$ and $\mu=(\mu_1, \ldots, \mu_n)$ be partitions of at most
$n$ parts for a positive integer $n$ such that the
skew partition $\lambda/\mu$ is a horizontal strip; i.e. $\lambda_1
\geq \mu_1 \geq\lambda_2 \geq \mu_2 \geq \ldots \lambda_n \geq
\mu_n \geq \lambda_{n+1} = \mu_{n+1} = 0$. Following~\cite{Coskun1}, define
\begin{multline}
\label{definitionHfactor}
H_{\lambda/\mu}(q,p,t,b) \\
:= \prod_{1\leq i < j\leq n} 
\left\{\dfrac{(q^{\mu_i-\mu_{j-1}}t^{j-i})_{\mu_{j-1}-\lambda_j}
(q^{\lambda_i+\lambda_j}t^{3-j-i}b)_{\mu_{j-1}-\lambda_j}}
{(q^{\mu_i-\mu_{j-1}+1}t^{j-i-1})_{\mu_{j-1}-\lambda_j}(q^{\lambda_i
    +\lambda_j+1}t^{2-j-i}b)_{\mu_{j-1}-\lambda_j}}\right.\\
\left.\cdot 
\dfrac{(q^{\lambda_i-\mu_{j-1}+1}t^{j-i-1})_{\mu_{j-1}-\lambda_j}}
{(q^{\lambda_i-\mu_{j-1}}t^{j-i})_{\mu_{j-1}-\lambda_j}}\right\} 
\cdot\prod_{1\leq i <(j-1)\leq n} \!\!\!
\dfrac{(q^{\mu_i+\lambda_j+1}t^{1-j-i}b)_{\mu_{j-1}-\lambda_j}}
{(q^{\mu_i+\lambda_j}t^{2-j-i}b)_{\mu_{j-1}-\lambda_j}}
\end{multline}
and 
\begin{multline}
\label{definitionSkewW}
W_{\lambda/\mu}(x; q,p,t,a,b)
:= H_{\lambda/\mu}(q,p,t,b)\cdot\dfrac{(x^{-1}, ax)_\lambda
  (qbx/t, qb/(axt))_\mu}
{(x^{-1}, ax)_\mu (qbx, qb/(ax))_\lambda}\\
\cdot\prod_{i=1}^n\left\{\dfrac{\theta(bt^{1-2i}q^{2\mu_i})}{\theta(bt^{1-2i})}
  \dfrac{(bt^{1-2i})_{\mu_i+\lambda_{i+1}}}
{(bqt^{-2i})_{\mu_i+\lambda_{i+1}}}\cdot
t^{i(\mu_i-\lambda_{i+1})}\right\}
\end{multline}
where $q,p,t,x,a,b\in\mathbb{C}$. 
Note that $W_{\lambda/\mu}(x; q,p, t, a,b)$ vanishes unless $\lambda/\mu$ is a horizontal strip. 
The function
$W_{\lambda/\mu}(y, z_1, \ldots, 
z_\ell; q,p,t,a,b)$ is extended to $\ell+1$ variables $y, z_1, \ldots, z_\ell
\in\mathbb{C}$  
through the following recursion formula
\begin{multline}
\label{eqWrecurrence}
W_{\lambda/\mu}(y,z_1,z_2,\ldots,z_\ell;q, p, t, a, b) \\
= \sum_{\mu\prec \nu\prec \lambda} W_{\lambda/\nu}(yt^{-\ell};q, p, t, at^{2\ell},
bt^\ell) \, W_{\nu/\mu}(z_1,\ldots, z_\ell;q, p, t, a, b).
\end{multline}

\subsection{The Limiting $w_{\lambda/\mu}$ Function}

The Macdonald functions $W_\lambda$ are essentially equivalent to $BC_n$ abelian functions constructed independently in~\cite{Rains1}. 
The limiting cases defined above are closely related to the Macdonald polynomials~\cite{Macdonald1}, and interpolation Macdonald polynomials~\cite{Okounkov1}.

The following limiting 
 the basic (the $p=0$ case of the elliptic) $W$ functions
%$W_{\lambda/\mu}(x; q,t,a,b) = W_{\lambda/\mu}(x; q, 0,t,a,b)$ 
will be used in our constructions below. 
The existence of these limits can be seen from (the $p=0$ case of) the
definition~(\ref{definitionSkewW}),  the recursion
formula~(\ref{eqWrecurrence}) and the limit rule
\begin{equation}
\label{LimitRule}
  \lim_{a\rightarrow 0}\, a^{|\mu|} (x/a)_{\mu} 
= (-1)^{|\mu|}\, x^{|\mu|} t^{-n(\mu)} q^{n(\mu')} 
\end{equation}
where $\abs{\mu}=\sum_{i=1}^n \mu_i$ and $n(\mu) =
\sum_{i=1}^n (i-1) \mu_i$, 
and $n(\mu') =\sum_{i=1}^n \binom{\mu_i}{2}$. 
We denote $H_{\lambda/\mu}(q,t)=H_{\lambda/\mu}(q,0,t,0)$, and for $x\in \mathbb{C}$ define
\begin{multline}
\label{wdefn}
w_{\lambda/\mu}(x; q,t) := \lim_{s\rightarrow \infty} \left( s^{|\lambda|-|\mu|} \lim_{a\rightarrow 0} W_{\lambda/\mu}(x; q,t, a, as) \right) \\
= (-q/x)^{-|\lambda|+|\mu|} q^{ -n(\lambda') + n(\mu') }  H_{\lambda/\mu}(q, t) \dfrac{(x^{-1})_\lambda }{(x^{-1})_\mu }  
\end{multline}

The recurrence formula for $w_{\lambda/\mu}$ % $W^{\nearrow \hspace*{-6pt} s}_{\lambda/\mu}$ 
function turns out to be
\begin{equation}
\label{Wsuprec}
w_{\lambda/\mu}(y,z;q, t) \\
= \sum_{\nu\prec \lambda} t^{\ell(|\lambda|-|\nu|)}
w_{\lambda/\nu}(yt^{-\ell};q, t) \, 
w_{\nu/\mu}(z;q, t)
\end{equation}
Similarly, for $x\in \mathbb{C}$ define the dual function 
\begin{multline}
\hat w_{\lambda/\mu}(x; q,t) := \lim_{s\rightarrow 0} \left( \lim_{a\rightarrow 0} W_{\lambda/\mu}(x; q,t, a, as) \right) \\
= t^{-n(\lambda)+|\mu|+n(\mu) } H_{\lambda/\mu}(q, t) \dfrac{(x^{-1})_\lambda }
{(x^{-1})_\mu }   \hspace*{70pt}
\end{multline}
The recurrence formula for the dual $\hat w_{\lambda/\mu}(x; q,t) $
% $W^{\searrow \hspace*{-6pt} s}_{\lambda/\mu}$ function 
may be written as
\begin{equation}
\hat w_{\lambda/\mu}(y,z;q, t) \\
= \sum_{\mu\prec \nu\prec \lambda}
\hat w_{\lambda/\nu}(yt^{-\ell};q, t) \, 
\hat w_{\lambda/\mu}(z;q, t)
\end{equation}
for $y\in\mathbb{C}$ and $z\in\mathbb{C}^\ell$. We now recall some old, and derive some new basic properties of the $w$ function and its dual, and their connections.

\begin{cor}  Let $\mu$ be a partition of at most $n$ part, and $x=(x_1,\ldots, x_n)\in\mathbb{C}^n$. Then \\
\label{wdual}
\noindent
(1) The $w_\mu$ and its dual $\hat w_\mu$ are flipped versions of one another. That is,
\begin{equation}
\label{flipiden}
\hat w_\mu (x, q, t) = q^{-|\mu|} t^{-2 n(\mu) - (n - 1)|\mu| } w_\mu (1/x,1/q,1/ t) 
\end{equation}

\noindent
(2) The limit $\lim_{q\rightarrow 1} w_\mu(xt^{\delta(n)}; q, t) $ exists when denominators do not vanish. For the particular case  when $x=\lambda$ is a partition, we use the notation
\begin{equation}
\bar w_\mu(q^\lambda t^{\delta(n)}; 1 ,t) :=
\lim_{q\rightarrow 1} (1-q)^{-\mu_1}
 w_\mu(q^\lambda t^{\delta(n)}; q ,t)  
\end{equation}

\end{cor}

\begin{proof}
Both properties follow, by direct calculation, from  the definition~(\ref{wdefn}) of $w_{\lambda/\mu}$, the recurrence relation~(\ref{Wsuprec}) for $w_{\lambda/\mu}$, and limit formula~(\ref{LimitRule}), and the flip formula
\begin{equation}
\label{flip}
x^{|\mu|}   (x^{-1}, q, t)_\mu 
= (-1)^{|\mu|} q^{n(\mu')} t^{-n(\mu)} (x;q^{-1}, t^{-1})_{\mu}
\end{equation}
The proof also uses the result that, in the limit 
\begin{equation}
\lim_{q\rightarrow 1} H_{\lambda/\mu}(q,t) = \dfrac{(k)_m}{m!}
\end{equation}
where $k$ is the maximum of the list $k=\mathrm{max}\{ \lambda_1 - \lambda_2, \lambda_2-\lambda_3, \ldots, \lambda_{n-1}-\lambda_n \}$, 
and $m$ is the maximum of 
$m=\mathrm{max}\{ \lambda_1 - \mu_1, \mu_1-\lambda_2, \lambda_2-\mu_2, \ldots, \mu_{n-1}-\lambda_n \}$.
If the second list has a negative number (which means that $\lambda/\mu$ is not a horizontal strip), then $H_{\lambda/\mu}(q,t) =0$. 
\end{proof}

\begin{rem}
We will need the following properties from~\cite{CoskunDiscM} in what follows.\\

\noindent
(3)
If $z=xt^{\delta(n)}$, for some $x\in\mathbb{C}$, we get
\begin{multline}
w_{\mu}(xt^{\delta(n)};q,t)  %W^{s\uparrow}_{\mu}(xt^{\delta(n)};q,t) 
=  (-1)^{|\mu|} x^{|\mu|} t^{n(\mu)} q^{-|\mu|-n(\mu')} (x^{-1})_\mu
\! \prod_{1\leq i < j\leq n} \dfrac{(t^{j-i+1})_{\mu_i
-\mu_j} } {(t^{j-i})_{\mu_i -\mu_j} }  
\end{multline}
which, after flipping $q$ and $t$ and using the flip rule~(\ref{flip}), may be written as
\begin{equation}
\label{wsxtdelta}
w_{\mu}(xt^{\delta(n)};q,t) 
=  q^{-|\mu|}  (x;q^{-1}, t^{-1})_{\mu}
\! \prod_{1\leq i < j\leq n} \dfrac{(t^{j-i+1})_{\mu_i
-\mu_j} } {(t^{j-i})_{\mu_i -\mu_j} } 
\end{equation}

\noindent
(4)
The vanishing property of $W$ functions implies that%~\cite{CoskunDiscM} that 
\begin{equation}
\label{wvanish}
  w_{\mu}(q^\lambda t^\delta;q, t)=0
\end{equation}
when $\mu\not\subseteq \lambda$, where $\subseteq$ denotes the partial inclusion ordering.\\

\noindent
(5)
Let $\lambda$ be an $n$-part partition with $\lambda_n 
\neq 0$ and $0\leq k\leq \lambda_n$ for some integer $k$, and let 
$z=(x_1,\ldots,x_n)\in \mathbb{C}^n$. It was shown in~\cite{CoskunDiscM} that
\begin{equation}
\label{wrect}
w_{\bar k}(z;q,t) 
= q^{-nk} \prod_{i=1}^n (q^{1-k} x_i)_k 
\end{equation}
where $\bar k = (k,\ldots, k)\in \mathbb{C}^n$. \\

\noindent
(6)
With the same notation as above, we also have
\begin{equation}
\label{wnormal}
w_\mu(q^\mu t^{\delta(n)}; q, t) 
= q^{-|\mu|}\, t^{(n-1)|\mu|-2n(\mu)} \, ( qt^{n-1} )_{\mu} 
 \!\!\prod_{1\leq i < j\leq n} \!\!
\frac{(qt^{j-i-1})_{\mu_i-\mu_j} } {(qt^{j-i})_{\mu_i-\mu_j} }
\end{equation}

\end{rem}

\subsection{The Multiple $qt$-Binomial Coefficients}
\label{section1}
The multiple Stirling numbers we develop in this paper are  closely connected with binomial coefficients just as in the one dimensional case. Recall that the multiple $qt$-binomial coefficient is defined~\cite{CoskunDiscM} as
\begin{defn}
\label{qtbinomcoeffExt}
Let $z=(x_1,\ldots, x_n)\in\mathbb{C}^n$ and $\mu$ be an $n$-part partition. Then the $qt$-binomial coefficient is defined by  
\begin{equation}
\label{qtbinom}
\binom{z}{\mu}_{\!\!\!q,t} := \dfrac{ q^{|\mu|} t^{2n(\mu)+(1-n)|\mu| } }  { (qt^{n-1} )_\mu} \prod_{1\leq i<j \leq n} \left\{\dfrac{ (qt^{j-i})_{\mu_i-\mu_j} } {(qt^{j-i-1})_{\mu_i-\mu_j}  } \right\}  w_\mu(q^z t^{\delta(n)}; q, t)
\end{equation}
where $q,t\in\mathbb{C}$. It should be noted that this definition makes sense even for $\mu\in\mathbb{C}^n$. 
\end{defn}

Setting $t=q^\alpha$ and sending $q\rightarrow 1$ yields the multiple ordinary $\alpha$-binomial coefficients. 
For $n=1$, the definition reduces to that of the one dimensional $q$-binomial coefficients
\begin{equation}
\binom{n}{k}_{\!\!q} : = \dfrac{ (q)_n }{(q)_{n-k} (q)_k }
\label{eq:qBinom}
\end{equation}
which are also known as the Gaussian polynomials. These are studied extensively in the literature including but not limited to the works in~\cite{Andrews5, Andrews6, GasperR1, BerkovichW1, Gould2, Macdonald2, GarvanS1, ClarkI1}. 

Using the definition above, we write~\cite{CoskunDiscM} the terminating $qt$-binomial theorem in the form
\begin{equation}
\label{qt_binom_thmAlt}
 (x)_\lambda 
 = \sum_{\mu \subseteq \lambda} (-1)^{|\mu|} q^{n(\mu')} t^{-n(\mu)} \binom{\lambda}{\mu}_{\!\!\!q,t} x^{|\mu|} 
\end{equation}
This special case and its nonsymmetric analogues are also proved independently in~\cite{Okounkov1, Sahi1}, and studied in more recent works~\cite{Lascoux1}. 

For the special case $z=\bar x=(x,\ldots,x)\in\mathbb{C}^n$, we get
\begin{multline}
\label{binom_rect}
\binom{\bar x}{\mu}_{\!\!\!q,t} = 
\dfrac{ t^{2n(\mu)+(1-n)|\mu| } }  { (qt^{n-1} )_\mu} \,  (q^x;q^{-1}, t^{-1})_{\mu} \\
\cdot \prod_{1\leq i<j \leq n} \left\{\dfrac{ (qt^{j-i})_{\mu_i-\mu_j} } {(qt^{j-i-1})_{\mu_i-\mu_j}  }  \dfrac{(t^{j-i+1})_{\mu_i
-\mu_j} } {(t^{j-i})_{\mu_i -\mu_j} } \right\} \hspace{50pt}
% W^{\nearrow \hspace*{-6pt} s}_\mu(q^z t^{\delta(n)}; q, t) \\
\end{multline}
This is also established independently in~\cite{Lassalle1}. 

Another special case occurs when $\mu$ is a rectangular partition, that is $\mu=\bar k$. Using~(\ref{wrect}) we get
\begin{equation}
\label{binom_rect1}
\binom{z}{\bar k}_{\!\!\!q,t} =\prod_{i=1}^n \dfrac{ (q^{1-k} q^{x_i}t^{n-i} )_k  }  { (qt^{n-i} )_{k}}  
\end{equation}
In the particular case for $k=1$, the definition reduces to
\begin{equation}
\binom{z}{\bar 1}_{\!\!\!q,t} =\prod_{i=1}^n \dfrac{ ( q^{x_i}t^{n-i})_1  }  { (qt^{n-i} )_{1}}  
= \prod_{i=1}^n \dfrac{(1-q^{x_i} t^{n-i} )}{(1-qt^{n-i} )} 
\end{equation}

The last property we point out from~\cite{CoskunDiscM} is that 
\begin{equation}
\binom{\lambda}{\lambda}_{\!\!\!q,t} = \binom{\lambda}{\bar 0}_{\!\!\!q,t} = 1 
\end{equation}
where $\lambda$ is an $n$-part partition, and $\bar 0$ is the $n$-part partition whose parts are all zero. In addition, 
\begin{equation}
\label{binomvanish}
\binom{\lambda}{\mu}_{\!\!\!q,t} = 0 
\end{equation}
when $\mu\not\subseteq \lambda$ or $\bar 0\not\subseteq\mu$.

\subsection{The Multiple $qt$-Factorial Function}
We now recall another important extension~\cite{CoskunDiscM} that generalizes the one dimensional $q$-bracket and $q$-factorial polynomial to the multiple case as follows. 
\begin{defn}
Let $\mu$ be a partition of at most $n$ parts, $z=(x_1,\ldots,x_n)\in\mathbb{C}^n$ and $s\in\mathbb{C}^n$. Then 
\begin{multline}
\label{qtnumberShifted}
[z, s]_{\mu} = [z, s, n, q, t]_\mu \\
:= q^{|\mu|} 
\prod_{i=1}^n  \left\{ \dfrac{1}{(1-qt^{n-i} )^{\mu_i}} \right\}
\prod_{1\leq i<j \leq n} \left\{ \dfrac{ (t^{j-i})_{\mu_i-\mu_j} } {(t^{j-i+1})_{\mu_i-\mu_j}  } \right\}  \, w_\mu(s q^z t^{\delta(n)}; q, t)  % W^{\nearrow \hspace*{-6pt} s}_\mu(s q^z t^{\delta(n)}; q, t) \\
\end{multline}
 is called the $qt$-factorial (bracket) function. 
Note that the definition combines a multiplicative variable $s$, and an exponential variable $z$. Depending on the application we often set $z=\bar 0$ and write $\langle s\rangle_{\mu} = [\bar 0, s]_{\mu}$, or set $s=\bar 1$ and write $[ z]_{\mu} = [z, \bar 1]_{\mu}$. Using the identity~(\ref{wrect}) in the special case when $\mu= \bar 1$, we define the $qt$-bracket as
\begin{equation}
\label{qtnumber}
[z] =  [z, \bar 1,  n, q, t]_{\bar 1} 
% = \binom{z}{\bar 1}_{\!\!\!q,t} 
= \prod_{i=1}^n \dfrac{(1-q^{x_i} t^{n-i} )}{(1-qt^{n-i} )} 
\end{equation}
which is a multiple analogue of the classical $q$-number or $q$-bracket.
\end{defn}

The definition~(\ref{qtnumberShifted}) of $[ z]_{\mu}$ may also be written as 
\begin{multline}
\label{qtnumberShifted2}
[z]_{\mu} 
= t^{-2n(\mu)-(1-n)|\mu| } 
\prod_{i=1}^n  \left\{ \dfrac{ (qt^{n-i} )_{\mu_i} }{(1-qt^{n-i} )^{\mu_i}} \right\} \\
\cdot   \prod_{1\leq i<j \leq n} \left\{\dfrac{ (t^{j-i})_{\mu_i-\mu_j} } {(t^{j-i+1})_{\mu_i-\mu_j}  }     \dfrac{(qt^{j-i-1})_{\mu_i-\mu_j}  } { (qt^{j-i})_{\mu_i-\mu_j} } \right\} \binom{z}{\mu}_{\!\!\!q,t} 
\end{multline}
which reduces to~(\ref{qtnumber}) when $\mu=\bar 1$ by the identity~(\ref{binom_rect1}).

\begin{rem} 
The $qt$-factorial function satisfies the following properties: \\

\noindent
(a) 
Let $z=(x,\ldots, x)=\bar x  \in \mathbb{C}^n$ for a single variable $x\in\mathbb{C}$, then the $qt$-factorial function $[\bar x ]_{\mu}$ may be written as 
\begin{equation}
\label{spec_bracket}
[\bar x]_{\mu} 
=\prod_{i=1}^n  \left\{ \dfrac{1 }{(1-qt^{n-i} )^{\mu_i}} \right\} 
   (q^x;q^{-1}, t^{-1})_{\mu}
= \prod_{i=1}^n  \left\{ \dfrac{ (q^x t^{i-1};q^{-1})_{\mu_i} }{(1-qt^{n-i} )^{\mu_i}} \right\} 
\end{equation}
This definition reduces to the classical $q$-bracket in the one variable case. \\

\noindent
(b) Note that $(x;1/q,1/t)_\mu $, with the reciprocals of $q$ and $t$, corresponds to a multiple basic $qt$-analogue of the falling factorial $x_{\underline{n}}:= x (x-1) \cdots (x-(n-1))$ as opposed to the rising factorial or the Pochhammer symbol.  \\

\noindent
(c)
Setting $z=\mu$, and substituting the evaluation~(\ref{wnormal}) in~(\ref{qtnumberShifted}) above gives
\begin{multline}
[\mu]_{\mu} 
= t^{-2n(\mu)-(1-n)|\mu| } \,
\prod_{i=1}^n  \left\{ \dfrac{ (qt^{n-i} )_{\mu_i} }{(1-qt^{n-i} )^{\mu_i}} \right\} \\
\cdot   \prod_{1\leq i<j \leq n} \left\{\dfrac{ (t^{j-i})_{\mu_i-\mu_j} } {(t^{j-i+1})_{\mu_i-\mu_j}  }     \dfrac{(qt^{j-i-1})_{\mu_i-\mu_j}  } { (qt^{j-i})_{\mu_i-\mu_j} } \right\} \hspace{20pt}
\end{multline}
which is precisely the front factor in~(\ref{qtnumberShifted2}) as expected. Similar to the classical case, we may use the notation $\mu!= [\mu]_{\mu} $ and write
\begin{equation}
\label{qtnumberShifted3}
[z]_{\mu} 
= \mu! \, \binom{z}{\mu}_{\!\!\!q,t}  
\end{equation}
Note that in the particular case when $\mu=\bar k$ is a rectangular partition, the $\mu!$ reduces to a product of one dimensional quotients for each part.
\begin{equation}
\label{rect!}
\bar k !
= \prod_{i=1}^n  \left\{ \dfrac{ (qt^{n-i} )_{k} }{(1-qt^{n-i} )^{k}} \right\} 
\end{equation}

\noindent
(d) 
Recall that a multiple analogue $ E_{q,t}(z)$ of the exponential function $e^z$ was defined~\cite{Coskun1} by  
\begin{multline}
 E_{q,t}(z):=(-z)_{\infty^n}
 = \sum_{\mu\in P_n} \dfrac{z^{|\mu|} q^{n(\mu')} t^{n(\mu)+(1-n)|\mu| } }  { (qt^{n-1} )_\mu} \\ \cdot \prod_{1\leq i<j \leq n} \left\{\dfrac{ (qt^{j-i})_{\mu_i-\mu_j} } {(qt^{j-i-1})_{\mu_i-\mu_j}  }
  \dfrac{(t^{j-i+1})_{\mu_i -\mu_j} } {(t^{j-i})_{\mu_i -\mu_j} } \right\} 
\end{multline}
where $z\in \mathbb{C}$, and $P_n$ denotes the set of all partitions of length at most $n$, and $(z)_{\infty^n} = \prod_{i=1}^{n} (zt^{1-i})_{\infty}$. Using the $\mu!$ notation and setting $z=1$ gives a multiple analogue of the number $e$:
\begin{equation}
 E_{q,t}(1):=(-1)_{\infty^n}
 = \sum_{\mu\in P_n} \dfrac{ q^{n(\mu')} t^{-n(\mu)} }{\prod_{i=1}^n    (1-qt^{n-i} )^{\mu_i} }
\cdot \dfrac{1}{\mu !}
\end{equation}
similar to the classical case
\[  e=\sum_{m=0}^\infty \dfrac{1}{m!} \cdot \]

\end{rem}

\section{Recurrence Relations for Multiple Factorial Function}
We treat the recurrence relations for the multiple factorial function in a separate section here. We derive three distinct recurrences for the multiple binomial coefficients, and for multiple factorials $[z]_\lambda$: A recurrence relation with respect to the variable $z$, another recurrence with respect to the index $\lambda$, and a third recurrence with respect to the weight $|\lambda|$ of the index partition. 

\subsection{The variable $z$}
We can write the $W$--Jackson sum~\cite{Coskun1} in the form
\begin{multline}
\label{WJackson}
W_{\lambda}(z;q,t,at^{-2n},bt^{-n})\\
=\dfrac{(s)_{\lambda}(as^{-1}t^{-n-1})_{\lambda}}
{(qbs^{-1}t^{-1})_{\lambda}(qbt^ns/a)_{\lambda}}
\cdot \prod_{1\leq i < j\leq n}
\left\{\dfrac{(t^{j-i+1})_{\lambda_i-\lambda_j}
(qbt^{-i-j+1})_{\lambda_i+\lambda_j}}
{(t^{j-i})_{\lambda_i-\lambda_j} (qbt^{-i-j})_{\lambda_i +
    \lambda_j}} \right\}  \\
\cdot\sum_{\mu \subseteq \lambda}
  \dfrac{(bs^{-1}t^{-n})_{\mu} (qbt^n/a)_\mu}{(qt^{n-1})_\mu
(as^{-1}t^{-n-1})_\mu} \cdot \prod_{i=1}^n \left\{
  \dfrac{(1-bs^{-1}t^{1-2i}q^{2\mu_i})}{(1-bs^{-1}t^{1-2i})}(qt^{2i-2})^{\mu_i}
\right\} \\ \cdot \prod_{1\leq i < j \leq n}
\left\{ \dfrac{(t^{j-i})_{\mu_i -\mu_j} (qt^{j-i})_{\mu_i
      -\mu_j}}{(qt^{j-i-1})_{\mu_i -\mu_j}(t^{j-i+1})_{\mu_i
      -\mu_j}} \dfrac{(bs^{-1}qt^{-i-j})_{\mu_i+\mu_j}
    (bs^{-1}t^{-i-j+2})_{\mu_i+\mu_j}} {(bs^{-1}t^{-i-j+1})_{\mu_i+\mu_j}
(qbs^{-1}t^{-i-j+1})_{\mu_i+\mu_j}} \right\} \\
\cdot W_\mu(q^\lambda t^{\delta(n)};q,t,bt^{1-2n},bs^{-1}t^{-n})\cdot
W_\mu(zs;q,t,as^{-2}t^{-2n},bs^{-1}t^{-n})
\end{multline}
where $z\in\mathbb{C}^n$. Set $b=ar$ in this identity and send $a\rightarrow 0$, multiply both sides by $(rt^n)^{|\lambda|}$ and send $r\rightarrow \infty$ using the limit rule~(\ref{LimitRule}), and  to get
\begin{multline}
\label{wJackson}
w_{\lambda}(z;q, t) 
 = (-1)^{|\lambda|}\, q^{-|\lambda|-n(\lambda')} t^{n(\lambda)} s^{-|\lambda|} (s)_{\lambda} 
 \cdot \prod_{1\leq i < j\leq n}
\left\{\dfrac{(t^{j-i+1})_{\lambda_i-\lambda_j} }
{(t^{j-i})_{\lambda_i-\lambda_j}  } \right\}  \\
\cdot\sum_{\mu \subseteq \lambda}
  \dfrac{ (-1)^{|\mu|}\, q^{2|\mu|} t^{n(\mu)} q^{n(\mu')} }{(qt^{n-1})_\mu
} \cdot  \prod_{1\leq i < j \leq n}
\left\{ \dfrac{(t^{j-i})_{\mu_i -\mu_j} (qt^{j-i})_{\mu_i
      -\mu_j}}{(qt^{j-i-1})_{\mu_i -\mu_j}(t^{j-i+1})_{\mu_i
      -\mu_j}} \right\} \\
%\cdot \left( \lim_{r\rightarrow \infty}\lim_{a\rightarrow 0} W_\mu(q^\lambda t^{\delta(n)};q,t,art^{1-2n},ars^{-1}t^{-n}) \right) \cdot
\cdot \left( \lim_{a\rightarrow 0} W_\mu(q^\lambda t^{\delta(n)};q,t, a t^{1-2n},  as^{-1}t^{-n}) \right) \cdot
w_\mu(zs;q,t) s^{-|\mu|}   
\end{multline}
Shifting $z$ by $\,s^{-1}z$ in~(\ref{wJackson}) above, setting $s^{-1}=q^\alpha$ and  $z=q^zt^{\delta(n)}$ and using the flip rule
\begin{equation}
\label{flip1}
x^{-|\mu|}   (x; q, t)_{\mu}= 
(-1)^{|\mu|} q^{n(\mu')} t^{-n(\mu)} (x^{-1}; q^{-1}, t^{-1})_\mu 
\end{equation}
to simplify the front factors, and substituting the definition~(\ref{qtbinom}) into~(\ref{wJackson}) gives
\begin{multline}
\label{zrec1}
\binom{z+\alpha}{\lambda}_{\!\!\!q,t} 
% \label{wJackson}
% w_{\lambda}( q^{z+\alpha}  t^{\delta(n)} ;q, t) 
 = \prod_{1\leq i<j \leq n} \!\! \left\{\dfrac{ (qt^{j-i})_{\lambda_i-\lambda_j} }{(qt^{j-i-1})_{\lambda_i-\lambda_j}  } 
\dfrac{(t^{j-i+1})_{\lambda_i-\lambda_j} }{(t^{j-i})_{\lambda_i-\lambda_j}  } \right\}  \\
\cdot \dfrac{t^{2n(\lambda)+(1-n)|\lambda| } \, (q^{\alpha}; q^{-1}, t^{-1})_\lambda  }{ (qt^{n-1} )_\lambda}  
\cdot \sum_{\mu \subseteq \lambda}
 (-1)^{|\mu|}\, q^{(\alpha+1)|\mu|+n(\mu')}  t^{-n(\mu)-(1-n)|\mu| }  \\
\cdot \prod_{1\leq i<j \leq n} \left\{\dfrac{ (t^{j-i})_{\mu_i-\mu_j} } {(t^{j-i+1})_{\mu_i-\mu_j}  }  \right\}  
\left( \lim_{a\rightarrow 0} 
W_\mu(q^\lambda t^{\delta(n)};q,t,at^{1-2n}, aq^{\alpha} t^{-n}) \right)     
\binom{z}{\mu}_{\!\!\!q,t}
\end{multline}
This identity defines a general recurrence relation for the multiple binomial coefficients with respect to the $z$ variable. Choosing $\alpha\in\mathbb{C}$ 
and $z\in\mathbb{C}^n$ properly yields interesting special cases. For an $n$ partition  $\lambda$, let $\lambda^i$ denote $ \lambda^i= \lambda+e^i$ for $e^i=(e^i_1, \ldots, e^i_j, \ldots, e^i_n)$ with $e^i_i=1$, and $e^i_j=0$ for $i\neq j$. Setting $z+\alpha = \lambda^i$ in~(\ref{zrec1}), for example, gives a recurrence for 
\[ \binom{\lambda^i}{\lambda}_{\!\!\!q,t} = \binom{(\lambda_1,\ldots, \lambda_{i-1},\lambda_i+1, \lambda_{i+1}, \ldots, \lambda_n) }{(\lambda_1,\ldots, \lambda_{i-1},\lambda_i, \lambda_{i+1}, \ldots, \lambda_n)}_{\!\!\!q,t} \] 
for any $1\leq i \leq n$.

If we further substitute~(\ref{qtnumberShifted2}) for the multiple binomial coefficients on both sides, we get a recurrence for the multiple bracket function again with respect to the $z$ variable as follows.
\begin{multline}
% \label{wJackson}
%w_{\lambda}( q^{z+\alpha}  t^{\delta(n)} ;q, t) 
\dfrac{\prod_{i=1}^n  \left\{  (1-qt^{n-i} )^{\lambda_i} \right\} }
{(q^{\alpha}; q^{-1}, t^{-1})_\lambda  } \, [z+\alpha]_{\lambda} \\
 =  \sum_{\mu \subseteq \lambda}
  \dfrac{ (-1)^{|\mu|}\, q^{(\alpha+1) |\mu|}  t^{n(\mu)} q^{n(\mu')} }{(qt^{n-1})_\mu
} \cdot  \prod_{1\leq i < j \leq n} \!\!\!
\left\{ \dfrac{ (qt^{j-i})_{\mu_i
      -\mu_j}}{(qt^{j-i-1})_{\mu_i -\mu_j}  } \right\} \\
\cdot \left(  \lim_{a\rightarrow 0} W_\mu(q^\lambda t^{\delta(n)};q,t,at^{1-2n},aq^{\alpha} t^{-n}) \right)     \,
  \prod_{i=1}^n  \left\{  (1-qt^{n-i} )^{\mu_i} \right\} 
 [z]_{\mu}
\end{multline}
We may again specialize $z$ and $\alpha$ to derive interesting special cases. 

\subsection{The index $\lambda$}
We have shown~\cite{CoskunG1} that for $0<k\leq \lambda_n$
\begin{multline}
\label{commonfactor}
W_{\lambda}(z;q,t,a,b) 
=\prod_{1\leq i<j\leq n}\frac{(qbt^{j-2i})_{2k}}{(qbt^{j-1-2i})_{2k}} \\
\cdot \prod_{i=1}^n\frac{(z_i^{-1})_k (az_i)_k}{(qbz_i)_k (qb/(az_i))_k}\,
W_{\lambda-\bar k}(zq^{-k}; aq^{2k}, bq^{2k}, t, q)
\end{multline}
where $\lambda-\bar k$ denotes $\lambda-k^n=(\lambda_1-k,\ldots,\lambda_n-k)$. 
We shift $\mu$ by $\bar k$, replace $b$ by $as$, send $a\rightarrow 0$, multiply both sides by $s^{|\mu|+kn}$, and then send $s\rightarrow \infty$ in the last identity, and use
the limit rule~(\ref{LimitRule}) and the definition~(\ref{wdefn}) of $w$ function to get 
\begin{equation}
w_{\mu+\bar k}(z;q,t) 
= (-1)^{nk}  q^{-nk - n \binom{k}{2} } \prod_{i=1}^n z_i^{k}  (z_i^{-1})_k  \cdot
w_{\mu}(zq^{-k}; t, q)  
\end{equation}

We now expand the right hand side of this identity using~(\ref{wJackson}). Replace $z$ by $\,s^{-1}z$, set $s=q^k$ and  $z=q^zt^{\delta(n)}$, and substitute the definition~(\ref{qtbinom}) here. Then using the flip formula~(\ref{flip1}) and the calculation that $t^{2n(\lambda+\bar k)+(1-n)|\lambda+\bar k| } = t^{2n(\lambda)+(1-n)|\lambda| } $, we write 
\begin{multline}
%w_{\lambda+\bar k}(q^z  t^{\delta(n)}; q,t)  
 \binom{z}{\lambda+\bar k}_{\!\!\!q,t}  
= t^{2n(\lambda)+(1-n)|\lambda| }  
\prod_{i=1}^n (q^{z_i} t^{(n-i)}; q^{-1})_k  \\
\cdot 
\dfrac{ (q^{-k}; q^{-1}, t^{-1})_{\lambda}  } 
{ (qt^{n-1} )_{\lambda+\bar k}   } 
\!\!\! \prod_{1\leq i < j\leq n}  \!\!\!
\left\{\dfrac{(t^{j-i+1})_{\lambda_i-\lambda_j} }
{(t^{j-i})_{\lambda_i-\lambda_j}  } \dfrac{ (qt^{j-i})_{\lambda_i-\lambda_j} }
{(qt^{j-i-1})_{\lambda_i-\lambda_j}  }  \right\}  \\
\cdot\sum_{\mu \subseteq \lambda}
 (-1)^{|\mu|}\, q^{(k+1)|\mu|+n(\mu')} t^{-n(\mu)-(1-n)|\mu| }  
\!\!\! \prod_{1\leq i < j \leq n}
\left\{ \dfrac{(t^{j-i})_{\mu_i -\mu_j}  }{ (t^{j-i+1})_{\mu_i - \mu_j}} \right\} \\
\cdot \left( \lim_{a\rightarrow 0} 
W_\mu(q^\lambda t^{\delta(n)};q,t,at^{1-2n}, aq^{-k}t^{-n}) \right)  
\binom{z}{\mu}_{\!\!\!q,t}  
\end{multline}
This is a general recurrence relation for the multiple binomial with respect to the index partition $\mu$. If we replace the binomials by the corresponding bracket functions using~(\ref{qtnumberShifted2}), we get 
\begin{multline}
[z]_{\lambda+\bar k} \,
= \prod_{i=1}^n  \left\{ \dfrac{ (q^{z_i} t^{(n-i)}; q^{-1})_k  } { (1-qt^{n-i} )^{k} } \right\}  
\cdot \dfrac{ (q^{-k}; q^{-1}, t^{-1})_{\lambda}  } {   \prod_{i=1}^n  (1-qt^{n-i} )^{\lambda_i} } \\
\cdot\sum_{\mu \subseteq \lambda}
 (-1)^{|\mu|}\, q^{(k+1)|\mu|+n(\mu')} t^{n(\mu) }  
\left( \lim_{a\rightarrow 0} 
W_\mu(q^\lambda t^{\delta(n)};q,t,at^{1-2n}, aq^{-k}t^{-n}) \right)  \\
\cdot 
\prod_{i=1}^n  \left\{ \dfrac{(1-qt^{n-i} )^{\mu_i}}{ (qt^{n-i} )_{\mu_i} } \right\} 
 \prod_{1\leq i<j \leq n} \left\{ \dfrac{ (qt^{j-i})_{\mu_i-\mu_j} }{(qt^{j-i-1})_{\mu_i-\mu_j}  } \right\} \, [z]_{\mu}  
\end{multline}
Specializing the variables $z$ and the integer $k$ gives interesting special cases as before. 

\subsection{The weight $|\lambda|$.} 
Recall that a recurrence relation for the $qt$-binomial coefficients based on the weights of the indexing partitions is written~\cite{Coskun1} as follows: For an $n$-part partition $\lambda$ and $k\leq |\lambda|$, we have
\begin{multline}
\sum_{\mu\vdash k} q^{n(\mu')} t^{-n(\mu)} \binom{\lambda^i}{\mu}_{\!\!\!qt} \\
 =
 \sum_{\tau \vdash k} q^{n(\tau')} t^{-n(\tau)} \binom{\lambda}{\tau}_{\!\!\!qt} + 
\sum_{\nu \vdash (\!k-1\!)} \!\!\! q^{n(\nu')+\lambda_i } t^{-n(\nu)+1-i} \binom{\lambda}{\nu}_{\!\!\!qt} \hspace{30pt}
\end{multline}
where $\mu\subseteq\lambda^i = (\lambda_1, \ldots, \lambda_{i}+1, \ldots, \lambda_n)$ when it is a partition, and  $\tau,\nu\subseteq\lambda$. Here, $\mu\vdash k$ denotes that $\mu$ is a partition of the integer $k\geq 1$. Using~(\ref{qtnumberShifted3}), we write a recurrence formula for the $qt$-factorial function $[\lambda]_{\mu}$ in the following form:
\begin{multline}
\sum_{\mu\vdash k} q^{n(\mu')} t^{-n(\mu)} \, \mu! \, [\lambda^i]_{\mu}   \\
 =  \sum_{\tau \vdash k} q^{n(\tau')} t^{-n(\tau)} \, \tau! \, [\lambda]_{\tau}    + 
\sum_{\nu \vdash (\!k-1\!)} \!\!\! q^{n(\nu')+\lambda_i } t^{-n(\nu)+1-i} \, \mu! \, [\lambda]_{\mu}   
\end{multline}

\section{Summation Formulas for Multiple Factorial Functions}
The \ci~\cite{CoskunG1} for $\omega$ functions
\begin{equation}
  \omega_{\lambda/\mu}((sr)^{-1},sr,as^2, bs) 
=\sum_\nu
  \omega_{\lambda/\nu}(s^{-1},s,as^2,bs) \,
  \omega_{\nu/\mu}(r^{-1};r,a,b)  
\end{equation}
may be written in explicit form as  
\begin{multline}
\dfrac{(rs)_{\lambda} (asr^{-1})_{\lambda}} {(qbr^{-1})_{\lambda}
    (qbr/a)_{\lambda}} \dfrac {(qb)_{\lambda}
    (qb/a)_{\lambda}} {(s)_{\lambda} (as)_{\lambda}}
    \dfrac{ (qb/as)_{\mu}}{ (asr^{-1})_{\mu}}  \dfrac{
    (ar^{-1})_{\mu}} {(qb/a)_{\mu}}  \\      
\cdot W_{\mu} (t, bst^{2-2n}, br^{-1}t^{1-n}; q^{\lambda}t^{\delta(n)};q) \\
=  \sum_\nu q^{|\nu|} t^{2n(\nu)} \dfrac{ (qb/as)_{\nu}}{ (as)_{\nu}} \dfrac{
    (bt^{1-n})_{\nu}}{ (qt^{n-1})_{\nu}} \dfrac{(r)_{\nu}
    (ar^{-1})_{\nu}} {(qbr^{-1})_{\nu} (qbr/a)_{\nu}}\\     
\cdot  \prod_{i=1}^{n}\left\{ \dfrac{(1-bt^{2-2i} q^{2\nu_i})}
    {(1-bt^{2-2i})}  \right\} 
  \prod_{1\leq i< j \leq n} \left\{ \dfrac{ (qt^{j-i})_{\nu_i - \nu_j}
    } { (qt^{j-i-1})_{\nu_i - \nu_j} } \dfrac{
    (bt^{3-i-j})_{\nu_i + \nu_j} } { (bt^{2-i-j})_{\nu_i + 
    \nu_j} } \right\} \\ 
\cdot W_{\nu} (t, bst^{2-2n}, bt^{1-n};
    q^{\lambda}t^{\delta(n)};q) \, W_{\mu} (t, bt^{2-2n}, br^{-1}t^{1-n};
    q^{\nu}t^{\delta(n)};q)   
\end{multline}
A specialization of this result gave the elegant identity for $qt$-binomial coefficients:
For $n$-part partitions $\nu$ and $\mu$, we have~\cite{CoskunDiscM}
\begin{equation}
\label{doublebinom}
t^{-n(\mu)} q^{n(\mu')} 
\dfrac{ ( -1 )_{\nu}} { (-1 )_{\mu} } \,\binom{\nu}{\mu}_{\!\!\!q,t} 
= \sum_{\mu \subseteq
    \lambda \subseteq \nu } 
    t^{-n(\lambda)} q^{n(\lambda')} 
    \binom{\nu}{\lambda}_{\!\!\!q,t}
    \binom{\lambda}{\mu}_{\!\!\!q,t}
\end{equation}
Setting $\mu=\bar 1$ here gives 
\begin{equation}
%\label{doublebinom}
t^{-n(\mu)} q^{n(\mu')} 
\dfrac{ ( -1 )_{\nu}} { (-1 )_{\mu} } \, [\nu]_{q,t}
= \sum_{\mu \subseteq
    \lambda \subseteq \nu } 
    t^{-n(\lambda)} q^{n(\lambda')}  
    \binom{\nu}{\lambda}_{\!\!\!q,t}  [\lambda]_{q,t}
\end{equation}
which is a multiple analogue of the identity 
\begin{equation}
2^{n-1} n = \sum_{k=1}^n  \binom{n}{k}  k
\end{equation}
for the $qt$-bracket function~(\ref{qtnumberShifted}). 

We like to also write a multiple analogue of the identity that expresses the sum of first $n$ integers as
\begin{equation}
\label{Gauss-sum}
\sum_{i=1}^n k = \binom{n+1}{2} = \dfrac{n(n+1)}{2} 
\end{equation}
for the multiple bracket. One dimensional $q$-analogues of this identity appeared in several papers~\cite{Garrett2, Warnaar3, Schlosser1} recently. For the multiple analogue, we first recall another special case of the \ci\ written~{\cite{Coskun1}} in the form 
\begin{multline}
\label{strongci_limit}
\dfrac{ ( zs )_{\nu}} {(s )_{\nu}} 
    \dfrac{( z )_{\mu}}{ (zs )_{\mu} } \, s^{|\mu|} \,
w_{\mu} (q^{\nu}t^{\delta(n)}; q, t) \\
=  \sum_{\mu \subseteq
\lambda \subseteq \nu} 
 q^{|\lambda|} t^{2n(\lambda)} 
\dfrac{(z )_{\lambda}}{(qt^{n-1})_{\lambda} } 
\!\! \prod_{1\leq i< j \leq n} \hspace*{-5pt} \left\{ \dfrac{
    (qt^{j-i})_{\lambda_i - \lambda_j} } { (qt^{j-i-1})_{\lambda_i - \lambda_j} }
      \right\} \\
%\cdot W^{ab}_{\lambda} (q^{\nu}t^{\delta(n)}; q, t, s^{-1}t^{n-1} ) 
\cdot \left( \lim_{a\rightarrow 0} W_\lambda(q^\nu t^{\delta(n)};q,t,at^{1-2n},as^{-1}t^{-n}) \right) w_{\mu} (q^{\lambda}t^{\delta(n)}; q, t) 
\end{multline}
It follows from the recurrence~(\ref{eqWrecurrence}) that the $W$ function in the sum simplifies to 
\begin{equation}
\label{Wablimit}
 \left( \lim_{a\rightarrow 0} W_\lambda(q^\nu t^{\delta(n)};q,t,at^{1-2n},aq^{-1}t^{1-2n}) \right) \\
= \prod_{1\leq i < j\leq n} \dfrac{(t^{j-i+1})_{\lambda_i
-\lambda_j} } {(t^{j-i})_{\lambda_i -\lambda_j} }
\end{equation}
in the special case when $s=qt^{n-1}$. Note that this specialization removed the $W$ function that terminated the sum from above. 
Using the definition of multiple bracket function~(\ref{qtnumberShifted}) on the right and multiple binomial~(\ref{qtbinom}) on the left, we may write 
\begin{multline}
%\label{strongci_limit}
t^{-2n(\mu)-2(1-n)|\mu| }  ( z )_{\mu}
\dfrac{ ( zqt^{n-1} )_{\nu}} {(qt^{n-1} )_{\nu}} 
    \dfrac{ (qt^{n-1} )_\mu }{ (zqt^{n-1} )_{\mu} } \, 
 \prod_{1\leq i<j \leq n} \left\{\dfrac{(qt^{j-i-1})_{\mu_i-\mu_j}  }{ (qt^{j-i})_{\mu_i-\mu_j} }  \right\}  
\binom{\nu}{\mu}_{\!\!\!q,t} \\
%w_{\mu} (q^{\nu}t^{\delta(n)}; q, t) \\
=  \sum_{\mu \subseteq
\lambda \subseteq \nu} \dfrac{(z )_{\lambda}}{(qt^{n-1})_{\lambda} } 
\, q^{|\lambda|} t^{2n(\lambda)} 
\prod_{1\leq i< j \leq n} \hspace*{-5pt} \left\{ \dfrac{
    (qt^{j-i})_{\lambda_i - \lambda_j} } { (qt^{j-i-1})_{\lambda_i - \lambda_j} }
 \dfrac{(t^{j-i+1})_{\lambda_i
-\lambda_j} } {(t^{j-i})_{\lambda_i -\lambda_j} }     \right\} \\
\cdot q^{-|\mu|}  \prod_{i=1}^n  \left\{  (1-qt^{n-i} )^{\mu_i} \right\} 
\prod_{1\leq i<j \leq n} \left\{ \dfrac{(t^{j-i+1})_{\mu_i-\mu_j}  }{ (t^{j-i})_{\mu_i-\mu_j} } \right\}  \,  [\lambda]_{\mu}
\end{multline}
Setting $z=qt^{n-1}$ further removes the factor that terminates the sum from below as well, but we still get a finite sum 
\begin{multline}
\label{strongci_limit2}
 q^{|\mu|} t^{-2n(\mu)-2(1-n)|\mu| } 
\dfrac{ ( q^2t^{2(n-1)} )_{\nu}} {(qt^{n-1} )_{\nu}} 
    \dfrac{ (qt^{n-1} )^2_\mu }{ (q^2t^{2(n-1)} )_{\mu} } 
  \prod_{i=1}^n  \left\{ \dfrac{1}{ (1-qt^{n-i} )^{\mu_i} } \right\}  \\
\cdot   \prod_{1\leq i<j \leq n} \left\{\dfrac{(qt^{j-i-1})_{\mu_i-\mu_j}  }
{ (qt^{j-i})_{\mu_i-\mu_j} }  
\dfrac{ (t^{j-i})_{\mu_i-\mu_j} }{(t^{j-i+1})_{\mu_i-\mu_j}  } \right\}   
\binom{\nu}{\mu}_{\!\!\!q,t}   \\ 
%w_{\mu} (q^{\nu}t^{\delta(n)}; q, t) \\
=  \sum_{\mu \subseteq
\lambda \subseteq \nu}  q^{|\lambda|} t^{2n(\lambda)} 
\prod_{1\leq i< j \leq n} \hspace*{-5pt} \left\{ \dfrac{
    (qt^{j-i})_{\lambda_i - \lambda_j} } { (qt^{j-i-1})_{\lambda_i - \lambda_j} }
 \dfrac{(t^{j-i+1})_{\lambda_i
-\lambda_j} } {(t^{j-i})_{\lambda_i -\lambda_j} }  \right\}   [\lambda]_{\mu}
\end{multline}
In particular for the special case $\mu=\bar 1$, we get a multiple analogue of~(\ref{Gauss-sum}) as
\begin{multline}
%\label{strongci_limit}
q^n t^{(n-1) n} 
\dfrac{ ( q^2t^{2(n-1)} )_{\nu}} {(qt^{n-1} )_{\nu}} 
 \prod_{i=1}^n 
    \dfrac{  (1-q^{\nu_i} t^{n-i} ) }{ (1-q^2t^{2n-1-i} ) }  \\
%\binom{\nu}{\bar 1}_{\!\!\!q,t} \\
%w_{\mu} (q^{\nu}t^{\delta(n)}; q, t) \\
=  \sum_{\bar 1 \subseteq
\lambda \subseteq \nu}  q^{|\lambda|} t^{2n(\lambda)} 
 \prod_{1\leq i< j \leq n} \hspace*{-5pt} \left\{ \dfrac{
    (qt^{j-i})_{\lambda_i - \lambda_j} } { (qt^{j-i-1})_{\lambda_i - \lambda_j} }
 \dfrac{(t^{j-i+1})_{\lambda_i
-\lambda_j} } {(t^{j-i})_{\lambda_i -\lambda_j} }     \right\}  \,  [\lambda] 
\end{multline}

Note that other specializations of $z$ in~(\ref{strongci_limit2}) above, such as $z=0$, gives different versions of this multiple identity. 
Using similar methods and the recurrences for multiple fuctorail functions, one may also write multiple analogues of the sums of powers of brackets as well.

\section{Multiple basic and ordinary $qt$-Stirling numbers}
\label{mulipleSpecial}
In this section we review the definition and fundamental properties of the multiple 
Stirling numbers of the first and second kind indexed by partitions 
~\cite{CoskunDiscM}. 
The classical Stirling numbers of the first kind are defined to be the coefficients of the power functions $x^k$ in the expansion of the falling factorial $x_{\underline{n}}=x (x-1) \ldots (x-n+1)$ written as
\begin{equation}
x_{\underline{n}} = n! \binom{x}{n}  = \sum_{k=0}^n s_1(n,k) \,x^k  .
\end{equation}
The $q$-analogue of these numbers are defined in~\cite{Carlitz1} and their properties are studied in~\cite{Milne2, Gould1, Kim1, Zeng1} and the works of others. 

First, we recall~\cite{CoskunDiscM} the definition of the multiple $qt$-Stirling numbers generalizing the one dimensional $q$-analogues.
\begin{defn}
For an $n$ part partition $\lambda$, the $qt$-Stirling numbers of the first kind $s_1(\lambda,\mu)$ are defined by 
\begin{equation}
\label{defnS1}
% [x ]_\lambda 
[\bar x ]_\lambda 
= \sum_{\mu \subseteq \lambda}  q^{-n(\lambda')} t^{2n(\mu)-(n-1)|\mu|} s_1(\lambda,\mu)  \,  \prod_{i=1}^{\mu_1} [\bar x^{i}] 
\end{equation}
where $x\in\mathbb{C}$, $\bar x=\{x,\ldots, x\}\in\mathbb{C}^n$, and $\bar x^{i}=\{x,\ldots, x, q, \ldots, q\}\in\mathbb{C}^n$ with $\mu'_i$ copies of $x$ for the dual partition  $\mu'$. 
That is,
\begin{equation}
\label{defnS2}
\prod_{i=1}^{\mu_1} [\bar x^{i}] 
= \prod_{i=1}^n \dfrac{(1-q^{x} t^{n-i} )^{\mu_i}}{(1-qt^{n-i} )^{\mu_i}} 
\end{equation} 
Likewise, 
the $qt$-Stirling numbers of the second kind $s_2(\lambda,\mu)$ are defined by 
\begin{equation}
 \prod_{i=1}^n \dfrac{(1-q^x t^{n-i} )^{\lambda_i}}{(1-qt^{n-i} )^{\lambda_i}}  = 
\prod_{i=1}^{\lambda_1} [\bar x^{i}] 
= \sum_{\mu \subseteq \lambda} q^{n(\mu')} t^{-2n(\nu)+(n-1)|\nu|} s_2(\lambda,\mu)  \, [\bar x]_\mu 
%[x]_\mu 
\end{equation}
\end{defn}

We now recall the explicit formula for the $qt$-Stirling numbers given in~\cite{CoskunDiscM}, generalizing the one dimensional $q$-analogues. We will refer to certain identities derived in the proof several times later on. In addition, there are some notational changes adopted in this paper. Thus, we include the main steps of the proof here as well. 
\begin{thm}
For $n$-part partitions $\nu$ and $\mu$, an explicit formula for the $qt$-Stirling numbers of first and second kind $s_1(\lambda,\mu)$ and $s_2(\lambda,\mu)$ are given by
\begin{multline}
\label{s1exp}
s_1(\nu,\mu)=s_1(\nu,\mu, q, t)=  
\dfrac{ q^{n(\nu')} t^{-2n(\mu)+2(n-1)|\mu|} }
{\prod_{i=1}^n (1-qt^{n-i} )^{\nu_i-\mu_i}} \, f_\mu(q,t)
\\  \cdot 
\sum_{\mu \subseteq \lambda \subseteq \nu} \!\!\!
  u(\nu, \lambda, q, t)\, t^{(1-n)|\lambda|}  
\bar w_\mu(q^{-\lambda} t^{-\delta(n)}; 1 , 1/t)  
\end{multline}
and
\begin{multline}
\label{s2exp}
s_2(\nu,\mu) = s_2(\nu,\mu,q,t)  :=
\dfrac{ q^{-n(\mu')} t^{2n(\nu)+(1-n)|\nu|} } 
{\prod_{i=1}^n (1-qt^{n-i} )^{\nu_i-\mu_i}}  \\ \cdot  
  \sum_{\mu \subseteq \lambda \subseteq \nu} \!\!
 (-1)^{|\lambda|}  t^{n(\lambda) }  
\bar w_\lambda(q^\nu t^{\delta(n)}; 1 ,t) \, f_\lambda(q,t) \,  v(\lambda,\mu,q,t) 
\end{multline}
where $\hat w_\mu$ and $\bar w_\mu$ are as defined above in Corollary~\ref{wdual}, and $f(\mu)$, $u(\lambda,\mu)$ and $v(\lambda,\mu)$ are given by
\begin{equation}
f(\mu,q,t):=  
\prod_{i=1}^{n-1} \dfrac{ (t)_{\mu_i -\mu_{i+1}} }{ (t^{n-i} )_{\mu_i}}
%\prod_{1\leq i < j \leq n, j\neq i+1}
\prod_{\substack{1\leq i < j \leq n \\  j\neq i+1}}
\left\{ \dfrac{ ( t^{j-i})_{\mu_i -\mu_j}}{( t^{j-i-1})_{\mu_i -\mu_j} } \right\} 
 \mu_n! \prod_{i=1}^{n-1} (\mu_i -\mu_{i+1})! 
\end{equation}
\begin{equation}
\label{ufact}
u(\lambda,\mu,q, t)
:= \dfrac{ q^{|\mu|} t^{2n(\mu)}  }{(qt^{n-1})_\mu
} \prod_{1\leq i < j \leq n}
\left\{ \dfrac{ (qt^{j-i})_{\mu_i
      -\mu_j}}{(qt^{j-i-1})_{\mu_i -\mu_j} } \right\} 
 \hat w_\mu(q^\lambda t^{\delta(n)};q,t) 
\end{equation}
and
\begin{multline}
\label{vfact}
v(\lambda,\mu,q,t)
:=  (-1)^{|\mu|} q^{n(\mu')} t^{-n(\mu) } \binom{\lambda}{\mu}_{\!\!\!q,t}  \\ 
= \dfrac{ (-1)^{|\mu|} q^{|\mu|+n(\mu')} t^{n(\mu)+(1-n)|\mu|} }  { (qt^{n-1} )_\mu}   
\prod_{1\leq i<j \leq n} \left\{\dfrac{ (qt^{j-i})_{\mu_i-\mu_j} } {(qt^{j-i-1})_{\mu_i-\mu_j}  } \right\}  w_\mu(q^\lambda t^{\delta(n)}; q, t)  
\end{multline}
\end{thm}

\begin{proof}
Simplify the identity~(\ref{wJackson}) for the case $z=xt^{\delta(n)}$ using~(\ref{wsxtdelta}), and send $s\rightarrow \infty$ to get
\begin{equation}
\label{changeofbasis1}
(x;q^{-1}, t^{-1})_\lambda
  = \sum_{\mu \subseteq \lambda}
  u(\lambda, \mu) \, x^{|\mu|}  
\end{equation}
with the definition of $u(\lambda, \mu)$ above.
Similarly, apply the shifts $a\rightarrow a s^2$, $b\rightarrow b s$ and $x\rightarrow x/s$ in~(\ref{WJackson})
at the beginning, and follow the same steps except send $s$ to 0 to get
\begin{equation}
\label{changeofbasis2}
 x^{|\lambda|} 
 =   \sum_{\mu \subseteq \lambda}
v(\lambda, \mu)\, (x;q^{-1}, t^{-1})_{\mu}
\end{equation}
where $v(\lambda, \mu)$ is defined as in the theorem. 
It is clear, by a change of basis argument, that 
\begin{equation}
\label{inv_uandv}
\delta_{\nu\lambda} = \sum_{\mu \subseteq \lambda \subseteq \nu}
  u(\nu, \lambda) \,  v(\lambda,\mu) 
\end{equation}
Note that the \lhs~(\ref{changeofbasis2}) does not depend on $q$ or $t$. Now, flip the parameters $q\rightarrow 1/q, t \rightarrow 1/t$, 
take limit $q\rightarrow 1$, and multiply and divide the summand by 
$\prod_{i=1}^n 1/(1-qt^{n-i} )^{\mu_i}$ to get
\begin{equation}
\label{changeofbasis2b}
x^{|\lambda|} 
= \sum_{\mu \subseteq \lambda} 
\prod_{i=1}^n (1-qt^{n-i} )^{\mu_i}  \lim_{q\rightarrow 1} v(\lambda, \mu, 1/q, 1/t)  
\prod_{i=1}^{\mu_1} \langle \bar x^{i} \rangle 
\end{equation}
where $\langle \bar x^{i} \rangle$ defined as in~(\ref{defnS2}). Multiply and divide the summand in~(\ref{changeofbasis1}) by $t^{(n-1)|\mu|}$, 
substitute~(\ref{changeofbasis2b}) into~(\ref{changeofbasis1}) for $(xt^{n-1})^{|\mu|}$, 
multiply both sides of this latter identity by $\prod_{i=1}^n 1/(1-qt^{n-i} )^{\nu_i} $  to get
\begin{multline}
\prod_{i=1}^n  \left\{ \dfrac{1 }{(1-qt^{n-i} )^{\nu_i}} \right\}  (x;q^{-1}, t^{-1})_\nu 
= \sum_{\mu \subseteq \nu} \left(
 \prod_{i=1}^n \dfrac{ (1-qt^{n-i} )^{\mu_i} }{(1-qt^{n-i} )^{\nu_i}}  \right. \\ \cdot \left.
\sum_{\mu \subseteq \lambda \subseteq \nu} 
  u(\nu,\lambda, q, t) \, t^{-(n-1)|\lambda|} \Big( \lim_{q\rightarrow 1} v(\lambda,\mu, 1/q, 1/t) \Big)
\right) \prod_{i=1}^{\mu_1} \langle \bar x^{i} \rangle 
\end{multline}
Multiplying and dividing the summand now by $ q^{n(\nu')} t^{-2n(\mu)+(n-1)|\mu|}$ gives
\begin{equation}
\label{angledefnS1}
\langle \bar x\rangle_\nu \\ 
= \sum_{\mu \subseteq \nu} q^{-n(\nu')} t^{2n(\mu)-(n-1)|\mu|} s_1(\nu,\mu) 
% \prod_{i=1}^{\mu_1} \langle \bar x^{i} \rangle 
\cdot \prod_{i=1}^n \dfrac{(1- x t^{n-i} )^{\mu_i}}{(1-qt^{n-i} )^{\mu_i}} 
\end{equation}
where $s_1(\nu,\mu)$ is as defined in the theorem. A similar sequence of calculations show that
\begin{multline}
 \prod_{i=1}^n \dfrac{(1- x t^{n-i} )^{\nu_i}}{(1-qt^{n-i} )^{\nu_i}} 
= \sum_{\mu \subseteq \nu} \left(
 \prod_{i=1}^n \dfrac{ (1-qt^{n-i} )^{\mu_i} }{(1-qt^{n-i} )^{\nu_i}}  \right. \\ \cdot \left.
\sum_{\mu \subseteq \lambda \subseteq \nu} 
\Big( \lim_{q\rightarrow 1} u(\nu, \lambda,1/q,1/t)  \Big) t^{(n-1)|\lambda|} \, v(\lambda,\mu,q,t) 
\right)  \langle x \rangle_\mu 
\end{multline}
Multiplying and dividing the summand now by $ q^{-n(\mu')} t^{2n(\nu)+(1-n)|\nu|}$ gives the explicit formula for the $qt$-Stirling numbers of the second kind
\begin{equation}
\label{angledefnS2}
 \prod_{i=1}^n \dfrac{(1- x t^{n-i} )^{\nu_i}}{(1-qt^{n-i} )^{\nu_i}} 
= \sum_{\mu \subseteq \nu} q^{n(\mu')} t^{-2n(\nu)+(n-1)|\nu|} s_2(\nu,\mu)  \cdot 
\langle \bar x \rangle_\mu 
\end{equation}
Finally, substituting $x\rightarrow q^x$ completes the proof. We conclude by simplifying the flips and limits that entered the formulas above. 

It follows immediately from~(\ref{flip}) that, if
\begin{equation}
h(\mu, q, t):=\prod_{1\leq i<j \leq n} \left\{\dfrac{ (qt^{j-i})_{\mu_i-\mu_j} } {(qt^{j-i-1})_{\mu_i-\mu_j}  } \right\} 
\;
\mathrm{and}
\;\;
g(\mu, q, t):= (qt^{n-1} )_\mu 
\end{equation}
then 
\begin{equation}
h(\mu, 1/q, 1/t)=t^{2 n(\mu) - (n - 1) |\mu|}  h(\mu, q, t)
\end{equation}
and
\begin{equation}
g(\mu, 1/q, 1/t)= (-1)^{|\mu|} q^{-|\mu| - n(\mu')} t^{n(\mu) - (n - 1) |\mu|} g(\mu, q, t)
\end{equation}
Thus, flipping the parameters give
\begin{multline}
u(\lambda,\mu,1/q, 1/t)
= \dfrac{ (-1)^{|\mu|} q^{|\mu| + n(\mu')}  t^{n(\mu) - (n - 1)|\mu| }  }
{ (qt^{n-1})_\mu}  \\ \cdot \prod_{1\leq i < j \leq n}
\left\{ \dfrac{ (qt^{j-i})_{\mu_i -\mu_j}}{(qt^{j-i-1})_{\mu_i -\mu_j} } \right\} 
w_\mu(q^\lambda t^{\delta(n)};q,t)
\end{multline}
and
\begin{multline}
v(\lambda,\mu,1/q,1/t) \\
= \dfrac{ t^{ (n - 1)|\mu| }  }  { (qt^{n-1} )_\mu}   
\prod_{1\leq i<j \leq n} \left\{\dfrac{ (qt^{j-i})_{\mu_i-\mu_j} } {(qt^{j-i-1})_{\mu_i-\mu_j}  } \right\} \cdot   w_\mu(q^{-\lambda} t^{-\delta(n)}; 1/q, 1/t) 
\end{multline}
Multiply and divide both by $(1-q)^{\mu_1}$, and pass the limit as $q\rightarrow 1$ to get
\begin{multline}
\lim_{q\rightarrow 1} u(\lambda,\mu,1/q, 1/t)  \\
= (-1)^{|\mu|}  t^{n(\mu) - (n - 1)|\mu| } 
\prod_{i=1}^{n-1} \dfrac{ (t)_{\mu_i -\mu_{i+1}} }{ (t^{n-i} )_{\mu_i}}
%\prod_{1\leq i < j \leq n, j\neq i+1}
\prod_{\substack{1\leq i < j \leq n \\  j\neq i+1}}
\left\{ \dfrac{ ( t^{j-i})_{\mu_i -\mu_j}}{( t^{j-i-1})_{\mu_i -\mu_j} } \right\}  \\
\cdot \lim_{q\rightarrow 1} \left(  (1-q)^{\mu_1} \dfrac{ 1}{ (q )_{\mu_n } } 
\prod_{i=1}^{n-1} \dfrac{1}{(q)_{\mu_i -\mu_{i+1}} }   \right)   
\bar w_\mu(q^\lambda t^{\delta(n)}; 1 ,t) 
\end{multline}
where $\bar w_\mu$ is defined as above. Using the limit rule~(\ref{LimitRule}), direct calculations give that
\begin{equation}
\lim_{q\rightarrow 1} \left(  (1-q)^{\mu_1} \dfrac{ 1}{ (q )_{\mu_n } } 
\prod_{i=1}^{n-1} \dfrac{1}{(q)_{\mu_i -\mu_{i+1}} }   \right)   
= \mu_n! \prod_{i=1}^{n-1} (\mu_i -\mu_{i+1})!
\end{equation}
Hence
\begin{equation}
\lim_{q\rightarrow 1} u(\lambda,\mu,1/q, 1/t)  
= (-1)^{|\mu|}  t^{n(\mu) - (n - 1)|\mu| }  
\bar w_\mu(q^\lambda t^{\delta(n)}; 1 ,t)  \, f(\mu,q,t)
\end{equation}
Similarly, 
\begin{equation}
\lim_{q\rightarrow 1} v(\lambda,\mu,1/q, 1/t)  
=  t^{ (n - 1)|\mu| }
\bar w_\mu(q^{-\lambda} t^{-\delta(n)}; 1 , 1/t)  \, f(\mu,q,t) 
\end{equation}
which completes the proof.
\end{proof}

Note that these explicit formulas allow us to extend the definition of $s_{1}(\nu,\mu)$ and $s_{2}(\nu,\mu)$ to any $\nu,\mu\in\mathbb{C}^n$. In other words, the indices do not have to be partitions. This is a new property even in the one dimensional case. 

\begin{rem}  
The immediate properties of the Stirling numbers established in~\cite{CoskunDiscM} are listed as follows: \\

\noindent
(a) The multiple Stirling numbers $s_1(\nu,\mu)$ and $s_2(\nu,\mu)$ admit explicit combinatorial formulas which are derived in the Theorem above. \\

\noindent
(b) %One Dimensional Case: 
These explicit formulas reduce to those for the $q$-Stirling numbers given by Kim  in~\cite{Kim1} for $n=1$. Moreover, sending $q\rightarrow 1$ in that case yields classical Stirling numbers of both types.\\ 

\noindent
(c) %Inversion: 
The matrix $m$ with entries $m_{\lambda\mu}=s_1(\lambda,\mu)$ is invertible in the sense of $V$ algebra defined at the beginning of Section \ref{back}, and its inverse is given by $m^{-1}_{\lambda\mu}=s_2(\lambda,\mu)$.  
More precisely, we have 
\begin{equation}
\label{inv}
\delta_{\nu\lambda} = \sum_{\mu \subseteq \lambda \subseteq \nu}
  s_1(\nu, \lambda) \,  s_2(\lambda,\mu) = \sum_{\mu \subseteq \lambda \subseteq \nu}
  s_2(\nu, \lambda) \,  s_1(\lambda,\mu) 
\end{equation}
which follows immediately from the inversion relation~(\ref{inv_uandv}). \\

\noindent
(d) % Diagonal Entries: 
Similar to the one dimensional case for the $q$-Stirling numbers, we have
\begin{equation}
\label{diagonal}
  s_1(\lambda, \lambda) = s_2(\lambda,\lambda) = 1
\end{equation}
for an arbitrary $n$-part partition $\lambda$. \\

\noindent
(e) %Additive Case: 
Setting $t=q^\alpha$ and sending $q\rightarrow 1$ gives the multiple ordinary $\alpha$-Stirling numbers of the first and second kind. 
\end{rem}

\section{Summation identities for multiple Stirling numbers}
We derive some additional new properties of the multiple Stirling numbers in this section. \\

\noindent
(1) 
First note that 
$\lim_{x\rightarrow 0} \langle x \rangle_\mu = \prod_{i=1}^n 1/(1-qt^{n-i} )^{\mu_i} $ follows readily from the formula~(\ref{spec_bracket}) written in the $qt$-angle bracket function $\langle \bar x\rangle_{\mu}$ as
\begin{equation}
\label{xbarBracket}
\langle \bar x\rangle_{\mu} 
=\prod_{i=1}^n  \left\{ \dfrac{1 }{(1-qt^{n-i} )^{\mu_i}} \right\} 
   (x;q^{-1}, t^{-1})_{\mu}
= \prod_{i=1}^n  \left\{ \dfrac{ (x t^{i-1};q^{-1})_{\mu_i} }{(1-qt^{n-i} )^{\mu_i}} \right\} 
\end{equation}
In the multiple case, setting $x=0$ in~(\ref{angledefnS1}) and~(\ref{angledefnS2}) respectively gives
\begin{equation}
 \prod_{i=1}^n \dfrac{1}{(1-qt^{n-i} )^{\nu_i}}  = \sum_{\mu \subseteq \nu} q^{-n(\nu')} t^{2n(\mu)+(n-1)|\mu|} s_1(\nu,\mu) 
 \prod_{i=1}^n \dfrac{1}{(1-qt^{n-i} )^{\mu_i}}
\end{equation}
and
\begin{equation}
 \prod_{i=1}^n \dfrac{1}{(1-qt^{n-i} )^{\nu_i}} 
= \sum_{\mu \subseteq \nu} q^{n(\mu')} t^{-2n(\nu)+(n-1)|\nu|} s_2(\nu,\mu)  \prod_{i=1}^n \dfrac{1}{(1-qt^{n-i} )^{\mu_i}}
\end{equation}
These appear to be new identities, even in the one dimensional case $n=1$. The identities may be interpreted as giving the eigenvectors of certain operators. \\

\noindent
(2) Note that, for an $n$-part partition $\nu$ with $\nu_n\neq 0$, the bracket function $\langle x\rangle_\nu$ has roots at $x=t^{1-j} q^{m_j}$ for $j=1,\ldots,n$, and $m_j=0,\ldots, \nu_j-1$. The limit bracket function $ \prod_{i=1}^n (1- xt^{n-i} )^{\mu_i}/(1-qt^{n-i} )^{\mu_i}$ has roots at $x=t^{1-j}$ for $j=1,\ldots,n$. Therefore, if
we set $x=t^{1-j}q^{m_j}$ (for some $m_j< \nu_{j}$) in~(\ref{angledefnS2}) we get
\begin{equation}
0 = \sum_{\mu \subsetneq \nu} q^{-n(\nu')} t^{2n(\mu)-(n-1)|\mu|} 
 \prod_{i=1}^n \dfrac{(1- q^{m_j} t^{1-j+n-i} )^{\mu_i}}{(1-qt^{n-i} )^{\mu_i}} \,
s_1(\nu,\mu) 
\end{equation}
where the summation is over all partitions $\mu\subsetneq \nu$, that is all partitions $\mu\subseteq \nu$ such that $\mu_{j}\leq m_j$. This inequality follows from the vanishing property of the $w$ functions~(\ref{wvanish}). 

In the particular case, setting $x=q$ in~(\ref{angledefnS1}) gives
\begin{equation}
0 = \sum_{\mu \subseteq \nu} q^{-n(\nu')} t^{2n(\mu)-(n-1)|\mu|} s_1(\lambda,\mu) 
\end{equation}
which is an analogue of $\sum_{k=0}^n s_1(n,k) =0$
in the classical case.

In another special case, setting $x=t^{1-j}$ (i.e., $m_j=0$) in~(\ref{angledefnS1}) and~(\ref{angledefnS2}) would amount to vanishing of all factorial functions $ \langle x\rangle_\mu $ except the ones in whose index the $j$-th part (thus all parts $j+1,\ldots, n$ after $j$) are 0. That is, the factorial functions will be nonzero only for partitions such as $\mu=(\mu_1,\ldots,\mu_{j-1},0,\ldots,0)$, and others will vanish. This is particularly interesting, for the substitution $x=t^{1-j}$ the limit factorial functions $ \prod_{i=1}^n (1- xt^{n-i} )^{\mu_i}/(1-qt^{n-i} )^{\mu_i}$ also vanish except  for partitions $\mu=(\mu_1,\ldots,\mu_{n-j},0,\ldots,0)$. 

Therefore, setting  $x=t^{1-j}$ in~(\ref{angledefnS1}) and~(\ref{angledefnS2})  respectively gives
\begin{equation}
\label{rootsS1}
%\langle t^{1-j} \rangle_\nu = 
0 = \sum_{\mu \subsetneq \nu} q^{-n(\nu')} t^{2n(\mu)-(n-1)|\mu|} s_1(\nu,\mu) 
% \prod_{i=1}^{\mu_1} \langle \bar x^{i} \rangle 
 \prod_{i=1}^n \dfrac{(1- t^{1-j+n-i} )^{\mu_i}}{(1-qt^{n-i} )^{\mu_i}} 
\end{equation}
where the sum is over all partitions $\mu\subsetneq \nu$ such that $\mu=(\mu_1,\ldots,\mu_{j-1},0,\ldots,0)$. Likewise, 
\begin{equation}
\label{rootsS2}
% \prod_{i=1}^n \dfrac{(1- t^{1-j+n-i} )^{\nu_i}}{(1-qt^{n-i} )^{\nu_i}} = 
0 = \sum_{\mu \subsetneq \nu} q^{n(\mu')} t^{-2n(\nu)+(n-1)|\nu|} s_2(\nu,\mu)  \, \langle t^{1-j} \rangle_\mu 
\end{equation}
where the sum is over all partitions $\mu\subsetneq \nu$ such that $\mu=(\mu_1,\ldots,\mu_{n-j},0,\ldots,0)$. 
The particular cases when $j=1$ in~(\ref{rootsS1}) and $j=n$ in~(\ref{rootsS2}) show that the multiple Stirling numbers vanish when $\mu = \bar 0 =(0,\ldots, 0) \in \mathbb{C}^n$ as in the classical case. That is, 
\[s_1(\lambda, \bar 0) = s_2(\lambda, \bar 0) =0\] 
for any $n$-part partition $\lambda$ with $\lambda_n\neq 0$. \\

\noindent
(3) Note that setting $x=t^{1-j} q^{m_j}$ for $j=1,\ldots,n$, and $m_j \geq \nu_j$ in~(\ref{angledefnS1}) and~(\ref{angledefnS2}) will not vanish the bracket functions. In particular, setting $\nu = \bar x =  \bar k = (k,\ldots, k)\in \mathbb{Z}^n$ in~(\ref{angledefnS1}) and using~(\ref{rect!}) gives % for $k>1$ gives
\begin{multline}
\bar k !  % = [\bar k ]_{\bar k} 
=\prod_{i=1}^n  \dfrac{ (qt^{n-i} )_{k} }{(1-qt^{n-i} )^{k}}    
= \sum_{\mu \subseteq \bar k}  q^{-n(\lambda')} t^{2n(\mu)-(n-1)|\mu|} s_1(\bar k,\mu)  \,  
\prod_{i=1}^n \dfrac{(1-q^k t^{n-i} )^{\mu_i}}{(1-qt^{n-i} )^{\mu_i}}  
\end{multline}
which is an analogue of the classical identity 
\[ k! = % \sum_{k=0}^n (-1)^{n-k} s_1(n,k) = 
\sum_{m=0}^k s_1(k,m) \,k^m \]
for the special case when $\nu = \bar k$ is a rectangular partition. \\

\noindent
(4) Recall that, 
$ s_1(n,m) = -s_2(n,m) = -\binom{n}{2}$  
when $n-m=1$ for the classical Stirling numbers. Similarly, if the index partitions satisfy $|\lambda|-|\tilde \lambda|=1$, we have that 
\[s_1(\lambda,\tilde \lambda) = - s_2(\lambda, \tilde \lambda) \]
exactly as in the one dimensional case. 

The proof follows easily from the inversion~(\ref{inv}) relation, and the observation that there are only two partitions between $\lambda$ and $\tilde \lambda$ under the inclusion ordering, namely the two partitions themselves. That is,
\begin{equation}
0= \delta_{\lambda \tilde \lambda} = \sum_{\tilde \lambda \subseteq \mu \subseteq \lambda}
  s_1(\lambda,\mu) \,  s_2(\mu,\tilde \lambda) 
\end{equation}
which implies that $ s_1(\lambda,\lambda) \,  s_2(\lambda,\tilde \lambda) =- s_1(\lambda,\tilde \lambda) \,  s_2(\tilde \lambda,\tilde \lambda)$. That the diagonal entries of both type of multiple $qt$-Stirling numbers are 1 by~(\ref{diagonal}) is now enough to conclude. \\

\section{The $qt$-Lah Numbers}

The classical Lah numbers are defined to be the connection coefficients in the expansion 
\begin{equation}
x^{\overline{n}} = \sum_{k=0}^n L(n,k) \, x_{\underline{k}}
\end{equation}
where $x_{\underline{n}}=x (x-1) \ldots (x-n+1)$ denotes the falling factorial as before, and $x^{\overline{n}}=x (x+1) \ldots (x+n-1)$ the rising factorial. 
Various $q$-analogues of these numbers are developed in one dimensional case in~\cite{GarsiaR1, Wagner2} and others. 

We now give the definition of multiple $qt$-Lah numbers in terms of the multiple factorial function and its flipped version.
\begin{defn}
Let $[\bar x]^\lambda$ denote the multiple analogue of the rising factorial. That is,
\begin{equation}
\label{rising}
[\bar x]^\lambda = [\bar x, q, t]^\lambda := [\bar x, q^{-1}, t^{-1}]_\lambda
\end{equation}
For an $n$ part partition $\lambda$, the $qt$-Lah numbers $L(\lambda,\mu) = L(\lambda,\mu, q, t) $ are defined by 
\begin{equation}
\label{Lahnumber}
[\bar x]^\lambda  
= \sum_{\mu \subseteq \lambda} (-1)^{|\mu|} q^{-|\mu|+  2n(\mu') }  t^{-n(\mu)}    
L(\lambda,\mu, q, t)  \, [\bar x]_\mu 
\end{equation}
where $x\in\mathbb{C}$, and $\bar x=\{x,\ldots, x\}\in\mathbb{C}^n$ as before. 
\end{defn}

In one dimensional case, the Lah numbers admit some explicit representataions. We show that the same is true for the multiple Lah numbers next.  
\begin{thm} Let $\nu$ and $\mu$ be partitions with at most $n$-parts. Then
\begin{multline}
\label{explicitLah}
L(\nu,\mu)  = (-1)^{-|\nu|+|\mu|} q^{-|\nu|+|\mu|- 2n(\mu')} t^{n(\nu) + n(\mu)} \\
\cdot
 \prod_{i=1}^n  \left\{ (1-qt^{n-i} )^{-\nu_i+\mu_i}  \right\} 
 \sum_{\mu \subseteq \lambda\subseteq \nu}
 u(\nu, \lambda, q^{-1}, t^{-1}) \, v(\lambda, \mu,q,t) 
\end{multline}
where $u$ and $v$ factors are as defined in~(\ref{ufact}) and~(\ref{vfact}) above. 
\end{thm}

\begin{proof}
Multiply and divide the right hand side of~(\ref{changeofbasis1}) by 
$\prod_{i=1}^n (1-qt^{n-i} )^{\lambda_i}  $, flip the parameters $q$ and $t$, and use the definitions~(\ref{xbarBracket}) and~(\ref{rising}) to get  
\begin{equation}
%\label{changeofbasis1}
\langle x;q, t \rangle^\lambda
  = \sum_{\mu \subseteq \lambda}
\prod_{i=1}^n  \left\{ \dfrac{1 }{(1-q^{-1}t^{-(n-i)} )^{\lambda_i}} \right\} u(\lambda, \mu, q^{-1}, t^{-1}) \, x^{|\mu|}  
\end{equation}

Similarly, multiplying and dividing the right hand side of~(\ref{changeofbasis2})  by 
the factor $\prod_{i=1}^n (1-qt^{n-i} )^{\mu_i}  $ gives
\begin{equation}
%\label{changeofbasis2}
 x^{|\lambda|} 
 =   \sum_{\mu \subseteq \lambda}  v(\lambda, \mu,q,t)\,
\prod_{i=1}^n  \left\{ (1-qt^{n-i} )^{\mu_i}  \right\}  \langle x;q, t\rangle_{\mu}
\end{equation}
Combine the two to get
\begin{multline}
\langle x;q, t \rangle^\lambda = 
\sum_{\mu \subseteq \nu} \sum_{\mu \subseteq \lambda\subseteq \nu}  
\prod_{i=1}^n  \left\{ \dfrac{1 }{(1-q^{-1}t^{-(n-i)} )^{\nu_i}} \right\} u(\nu, \lambda, q^{-1}, t^{-1})  \\ 
\cdot v(\lambda, \mu,q,t)\,
\prod_{i=1}^n  \left\{ (1-qt^{n-i} )^{\mu_i}  \right\}  \langle x;q, t\rangle_{\mu}
\end{multline}
Setting $x\rightarrow q^x$ and comparing the last identity to the definition of multiple Lah numbers~(\ref{Lahnumber}) gives the desired result. 
\end{proof}

\subsection{Properties of Lah Numbers}
We now establish a few fundamental properties of multiple Lah numbers starting with some special evaluations. \\

\noindent
(1) 
Set $\nu=\mu=\lambda$ in~(\ref{explicitLah}) to write 
\begin{multline}
%\label{explicitLah}
L(\lambda,\lambda)  = (-1)^{-|\lambda|+|\lambda|} q^{-|\lambda|+|\lambda|- 2n(\lambda')} t^{n(\lambda) + n(\lambda)} \\
\cdot
 \prod_{i=1}^n  \left\{ (1-qt^{n-i} )^{-\lambda_i+\lambda_i}  \right\} 
 u(\nu, \lambda, q^{-1}, t^{-1}) \, v(\lambda, \mu,q,t) 
\end{multline}
or
\begin{equation}
\label{specialLah1}
L(\lambda,\lambda)  =  q^{- 2n(\lambda')} t^{2n(\lambda) } 
 u(\lambda, \lambda, q^{-1}, t^{-1}) \, v(\lambda, \lambda,q,t) 
\end{equation}
It is clear from the definition~(\ref{vfact}) that
\begin{equation}
v(\lambda,\lambda,q,t)
=  (-1)^{|\lambda|} q^{n(\lambda')} t^{-n(\lambda) } 
\end{equation}
Likewise, the identities and show that  
\begin{multline}
\hspace*{15pt} \hat w_\lambda(q^\lambda t^{\delta(n)}; q, t) \\
=  (-1)^{|\lambda|}\, t^{-n(\lambda)} q^{-|\lambda|-n(\lambda')} 
(qt^{n-1} )_\lambda 
\!\! \prod_{1\leq i < j\leq n} \!\! \frac{(qt^{j-i-1})_{\lambda_i-\lambda_j} }
{(qt^{j-i})_{\lambda_i-\lambda_j} } \hspace*{15pt}
\end{multline}
Setting $\lambda=\mu$, and substituting the last evaluation in~(\ref{ufact}), and flipping the parameters $q$ and $t$ shows 
\begin{equation}
u(\lambda,\lambda, 1/q, 1/t)
= (-1)^{|\lambda|}\,  t^{-n(\lambda)}  q^{n(\lambda')} 
\end{equation}
Putting these into~(\ref{specialLah1}) now gives that $L(\lambda,\lambda)  = 1$,
as in the classical case. \\

\noindent
(2) 
In the view of property~(\ref{binomvanish}), the definition~(\ref{vfact}), and the explicit formula~(\ref{explicitLah}), we see that $L(\lambda,\mu) = 0$ when $\mu\not\subseteq \lambda$. 

\noindent
(3)
We now derive a multiple analogue of the closed formula 
\begin{equation}
L(n,m)  = \binom{n}{m} \dfrac{(n-1) !}{(m-1) !} %= \binom{n-1}{m-1} \dfrac{n !}{m !}
\end{equation}
We first write a generalization of the identity~(\ref{doublebinom}). 

\begin{lem}
Let $r\in\mathbb{C}$, and $\nu$ is an $n$ part partition. Then 
\begin{multline}
(-1)^{|\mu|}\, r^{|\mu|}  t^{-n(\mu)} q^{n(\mu')}  
\dfrac{ (r)_{\nu}} {(r )_{\mu} }  
w_{\mu} (q^{\nu}t^{\delta(n)}; q, t) \\
= \sum_{\substack{\lambda\\ \mu \subseteq
    \lambda \subseteq \nu }} \dfrac{ (-1)^{|\lambda|}\, r^{|\lambda|}  
    t^{n(\lambda)+(1-n)|\lambda|} q^{|\lambda|+n(\lambda')} } {(qt^{n-1})_{\lambda} } 
 \prod_{1\leq i< j \leq n} \left\{ \dfrac{ (qt^{j-i})_{\lambda_i - \lambda_j}
    } { (qt^{j-i-1})_{\lambda_i - \lambda_j} } \right\} \\
\cdot   w_{\lambda} (q^{\nu}t^{\delta(n)}; t, q) \,
    w_{\mu} (q^{\lambda}t^{\delta(n)}; t,q) 
\end{multline}
\end{lem}

\begin{proof}
Recall the following transformation identity for $w_\lambda$ functions~\cite{Coskun1}, which is obtained from a multiple analogue of Bailey's $_{10}\phi_9$ transformation formula from~\cite{CoskunG1}.
\begin{multline}
\dfrac{(s)_{\nu} } { (z)_{\mu} } 
(-1)^{|\mu|}\, z^{|\mu|}  t^{-n(\mu)} q^{n(\mu')} \\
\cdot \sum_{\substack{\lambda\\ \mu \subseteq
    \lambda \subseteq \nu }} \dfrac{
    (z)_{\lambda}} {  (qt^{n-1})_{\lambda} } \,      
 q^{|\lambda|} t^{2n(\lambda)}
  \prod_{1\leq i< j \leq n} \left\{ \dfrac{ (qt^{j-i})_{\lambda_i - \lambda_j}
    } { (qt^{j-i-1})_{\lambda_i - \lambda_j} } \right\} \hspace{30pt} \\ 
%    \cdot W^{ab}_{\lambda} (s^{-1}t^{n-1}; q^{\nu}t^{\delta(n)};t, q)\,
\cdot \left( \lim_{a\rightarrow 0} W_\lambda(q^\nu t^{\delta(n)};q,t,at^{1-2n},as^{-1}t^{-n}) \right) 
     w_{\mu} (q^{\lambda}t^{\delta(n)}; t, q)   \\
= \sum_{\substack{\lambda\\ \mu \subseteq
    \lambda \subseteq \nu }} \dfrac{ (-1)^{|\lambda|}\, z^{|\lambda|} s^{|\lambda|} 
    t^{n(\lambda)+(1-n)|\lambda|} q^{|\lambda|+n(\lambda')} } {(qt^{n-1})_{\lambda} } 
 \prod_{1\leq i< j \leq n} \left\{ \dfrac{ (qt^{j-i})_{\lambda_i - \lambda_j}
    } { (qt^{j-i-1})_{\lambda_i - \lambda_j} } \right\} \\
\cdot   w_{\lambda} (q^{\nu}t^{\delta(n)}; t, q) \,
    w_{\mu} (q^{\lambda}t^{\delta(n)}; t,q) 
\end{multline}
The sum on the left may be evaluated using the summation formula~(\ref{strongci_limit}) to get the result where we set $r=zs$. 
\end{proof}

This powerful result is more general than~(\ref{doublebinom}) because of the additional free variable $r$. If we set $r=t^{2n-2}$, we get
\begin{multline}
\label{newiden}
(-1)^{|\mu|}\,  t^{-n(\mu) }  q^{ n(\mu')}  
\dfrac{ (t^{2(n-1)})_{\nu}} {(t^{2(n-1)} )_{\mu} }  
\binom{\nu}{\mu}_{\!\!\!q,t}  \\
= \sum_{\substack{\lambda\\ \mu \subseteq
    \lambda \subseteq \nu }}  (-1)^{|\lambda|}\,  
    t^{-n(\lambda)+2(n-1)|\lambda|} q^{n(\lambda')}   
\binom{\nu}{\lambda}_{\!\!\!q,t} 
\binom{\lambda}{\mu}_{\!\!\!q,t} \hspace{30pt}
\end{multline}
which will be useful in writing our explicit formula that we give next.

\begin{thm}
With the notation as above, 
\begin{multline}
\label{closedLahexplicit}
L(\nu, \mu)  = (-1)^{|\nu|+|\mu|} q^{-|\nu|+|\mu|} t^{n(\nu) - n(\mu)} \\
\cdot \prod_{i=1}^n  \left\{ (1-qt^{n-i} )^{-\nu_i+\mu_i}  \right\}  
\dfrac{ (t^{2(n-1)})_{\nu}} {(t^{2(n-1)} )_{\mu} }  \,
\binom{\nu}{\mu}_{\!\!\!q,t}  \hspace{30pt}
\end{multline}
\end{thm}

\begin{proof}
Substituting~(\ref{ufact}) and~(\ref{vfact}) into the fomula~(\ref{explicitLah}) gives 
\begin{multline}
\label{Lah2}
L(\nu, \mu)  = (-1)^{|\nu| } q^{-|\nu|+|\mu|- n(\mu')} t^{n(\nu) } 
 \prod_{i=1}^n  \left\{ (1-qt^{n-i} )^{-\nu_i+\mu_i}  \right\}  \\
\cdot \sum_{\mu \subseteq \lambda\subseteq \nu}
\dfrac{ q^{-|\lambda|} t^{-2n(\lambda)}  }{(q^{-1}t^{-(n-1)}, q^{-1}, t^{-1})_\lambda
} \prod_{1\leq i < j \leq n}
\left\{ \dfrac{ (q^{-1}t^{-(j-i)})_{\lambda_i
      -\lambda_j}}{(q^{-1}t^{-(j-i-1)})_{\lambda_i -\lambda_j} } \right\} \\
\cdot \hat w_\lambda(q^{-\nu} t^{-\delta(n)};q^{-1},t^{-1})  \binom{\lambda}{\mu}_{\!\!\!q,t}  
\end{multline}
By definition~(\ref{flipiden}) we have 
\begin{equation}
%\label{flipiden}
\hat w_\lambda(q^{-\nu} t^{-\delta(n)};q^{-1},t^{-1})  
= q^{|\lambda|} t^{2 n(\lambda) + (n - 1)|\lambda| } 
w_\lambda (q^{\nu} t^{\delta(n)}; q, t) 
\end{equation}
and by the flip rule~(\ref{flip1}) we get
\begin{equation}
 (q^{-1}t^{-(n-1)}, q^{-1}, t^{-1})_\lambda  = (-1)^{|\lambda|} q^{-n(\lambda')-|\lambda|} t^{n(\lambda)+(1-n)|\lambda|}    (q t^{n-1}; q, t)_{\lambda}
\end{equation}
Observe also that
\begin{multline}
\prod_{1\leq i < j \leq n} \!\!
\left\{ \dfrac{ (q^{-1}t^{-(j-i)})_{\lambda_i
      -\lambda_j}}{(q^{-1}t^{-(j-i-1)})_{\lambda_i -\lambda_j} } \right\} 
= t^{2n(\lambda)+(1-n)|\lambda|}  \!\!\!
\prod_{1\leq i < j \leq n} \!\!
\left\{ \dfrac{ (qt^{j-i})_{\lambda_i
      -\lambda_j}}{(qt^{j-i-1})_{\lambda_i -\lambda_j} } \right\} 
\end{multline}
Therefore, putting the multiple binomial~(\ref{qtbinom}) into~(\ref{Lah2}) gives
\begin{multline}
L(\nu, \mu)  = (-1)^{|\nu| } q^{-|\nu|+|\mu|- n(\mu')} t^{n(\nu) } 
 \prod_{i=1}^n  \left\{ (1-qt^{n-i} )^{-\nu_i+\mu_i}  \right\}  \\
\cdot \sum_{\mu \subseteq \lambda\subseteq \nu}
(-1)^{|\lambda|} 
q^{n(\lambda')} t^{-n(\lambda) + 2(n - 1)|\lambda| }   
\binom{\nu}{\lambda}_{\!\!\!q,t}
  \binom{\lambda}{\mu}_{\!\!\!q,t}  
\end{multline}
Finally, we use~(\ref{newiden}) to evaluate the sum which gives the desired result. 
\end{proof}

\noindent
(4)
The special evaluations we proved earlier, such as $L(\lambda, \lambda)=1$ and $L(\lambda, \mu)=0$ for $\mu\not\subseteq \lambda$ can be verified easily using this closed formula. Other special evaluations are also possible. For example, setting $\mu=\bar 1$ gives 
\begin{multline}
L(\nu, \bar 1)  = (-1)^{|\nu|+n} q^{-|\nu|+n} t^{n(\nu) - \binom{n}{2}}  \\
\cdot  (t^{2(n-1)})_{\nu} \prod_{i=1}^n  \left\{ \dfrac{(1-qt^{n-i} )^{-\nu_i+1} }{(1-t^{2n-1-i} ) } 
\dfrac{(1-q^{\nu_i} t^{n-i} )}{(1-qt^{n-i} )}   \right\}  
\end{multline}
which is an analogue of 
\[ L(n,1) = \binom{n}{1} \dfrac{(n-1)!}{(1-1)! } = n! \]
This may be viewed as an alternative definition for $\nu !$. \\

\noindent
(5)
The final property we derive here is that the Lah numbers are self inverse of themselves in the sense that
\begin{equation}
x^{\overline{n}} = \sum_{k=0}^n L(n,k) \, x_{\underline{k}}, \quad\mathrm{and}\quad 
x_{\underline{n}} = \sum_{k=0}^n (-1)^{n-k} L(n,k) \, x^{\overline{k}} .
\end{equation}
Note that flipping the parameters $q$ and $t$ in the identity~(\ref{Lahnumber}) gives the inverse result
\begin{equation}
\label{invLahnumber}
[\bar x]_\lambda  
= \sum_{\mu \subseteq \lambda} (-1)^{|\mu|} q^{|\mu|-2n(\mu') }  t^{n(\mu)}    L(\lambda,\mu, t^{-1}, q^{-1})  \, [\bar x]^\mu 
\end{equation}
so that the matrices defined by the entries $ L(\lambda,\mu, t, q)$ and $ L(\lambda,\mu, t^{-1}, q^{-1}) $ respectively are inverses of each other in the sense of Section~\ref{back}. 

\section{Conclusion }
We have derived several interesting identities for the multiple $qt$-factorial functions, multiple $qt$-Stirling numbers, and multiple $qt$-Lah numbers in the present paper. 
We will construct additional properties such as other recurrence relations they satisfy, their explicit evaluations in various other special cases, their combinatorial interpretations, other generating functions they satisfy, and their connections to different families of multiple combinatorial numbers in an upcoming article. The Stirling and Lah numbers have interesting connections to various branches in mathematics such as the one expressed in the classical Dobinski's formula. Such relations will also be formulated in that paper.

\end{document}